\documentclass[final, 3p, fleqn, authoryear]{elsarticle} 
\usepackage{graphicx}[dvipdfmx] 
\makeatletter
\newcommand{\E}[1]{#1^{\mathrm{UE}}}
\renewcommand{\S}[1]{#1^{\mathrm{SO}}}

\newcommand{\RM}[1]{#1^{\mathrm{RM}}}
\newcommand{\PRM}[1]{#1^{\mathrm{PRM}}}

\newcommand{\PBP}[1]{#1^{\mathrm{PBP}}}

\newcommand{\RP}[1]{#1^{\mathrm{RP}}}
\newcommand{\PRP}[1]{#1^{\mathrm{PRP}}}

\newcommand{\bracket}[1]{\left\{ #1 \right\}}

\newcommand{\problemname}[1]{\mbox{[#1]}\hspace{\fill} \notag}

\makeatother

\usepackage{soul}
\usepackage{chngcntr}
\usepackage{apptools}
\AtAppendix{\counterwithin{lem}{section}}
\usepackage{amsmath,amssymb,bm}

\DeclareMathOperator*{\argmin}{arg\,min}
\usepackage{mathtools}
\usepackage[ISO]{diffcoeff}
\usepackage{amsthm}
\theoremstyle{definition}
\newtheorem{thm}{Theorem}[section]

\newtheorem{dfn}{Definition}[section]
\newtheorem{lem}{Lemma}[section]
\newtheorem{asm}{Assumption}[section]
\newtheorem{pro}{Proposition}[section]

\newtheorem*{prf}{Proof} 

\usepackage{comment}

\usepackage{refcount}
\ExplSyntaxOn
\NewDocumentCommand{\replace}{mmm}
 {
  \marian_replace:nnn {#1} {#2} {#3}
 }

\tl_new:N \l_marian_input_text_tl
\tl_new:N \l_marian_search_tl
\tl_new:N \l_marian_replace_tl

\cs_new_protected:Npn \marian_replace:nnn #1 #2 #3
 {
  \tl_set:Nn \l_marian_input_text_tl { #1 }
  \tl_set:Nn \l_marian_search_tl { #2 }
  \tl_set:Nn \l_marian_replace_tl { #3 }
  \regex_replace_all:nnN { \b\u{l_marian_search_tl}\b } { \u{l_marian_replace_tl} } \l_marian_input_text_tl
  \tl_use:N \l_marian_input_text_tl
 }
\ExplSyntaxOff

\usepackage[hang,small,bf]{caption}
\usepackage[subrefformat=parens]{subcaption}
\usepackage{float}
\restylefloat{table}
\usepackage{multirow}
\usepackage{tabularx}

\usepackage{hyperref}
\hypersetup{colorlinks=true,urlcolor=blue}

\usepackage{cleveref}

\crefname{equation}{Eq.}{Eqs.}
\crefname{figure}{Figure}{Figures}
\crefname{subfigure}{Figure}{Figures}
\crefname{table}{Table}{Tables}

\crefname{thm}{Theorem}{Theorems}
\crefname{cor}{Corollary}{Corollary}
\crefname{dfn}{Definition}{Definition}
\crefname{lem}{Lemma}{Lemmas}
\crefname{asm}{Assumption}{Assumption}
\crefname{pro}{Proposition}{Proposition}
\crefname{cnj}{Conjecture}{Conjecture}

\crefname{section}{Section}{Sections}
\crefname{subsection}{Section}{Sections}

\newcommand{\ClK}{\mathcal{K}}

\newcommand{\ClN}{\mathcal{N}}

\newcommand{\ClT}{\mathcal{T}}

\newcommand{\Vt}[1]{\bm{#1}}

\newcommand{\Vtp}{\bm{p}}
\newcommand{\Vtq}{\bm{q}}

\newcommand{\Vtw}{\bm{w}}

\newcommand{\Vtrho}{\bm{\rho}}

\usepackage{txfonts}
\usepackage{pxfonts}

\usepackage{color}

\usepackage{tikz}
\usepackage{pgfplots}
\usetikzlibrary{shapes.geometric}
\usetikzlibrary{arrows.meta}
\usetikzlibrary{overlay-beamer-styles}
\usetikzlibrary{calc,3d}
\usetikzlibrary{matrix, positioning}
\usetikzlibrary{decorations.pathreplacing}

\tikzset{ 
	c1/.style = {trapezium, 
	trapezium angle = 60,
	minimum height=0.433cm, 
	minimum width=0.5cm, 
	trapezium stretches=true, 
	inner xsep=0.0pt, inner ysep=0.433pt, outer sep=0pt,
	fill=gray},
	c2/.style = {trapezium, 
	trapezium angle = 60,
	minimum height=0.433cm, 
	minimum width=1.0cm, 
	trapezium stretches=true, 
	inner xsep=0.5pt, inner ysep=0.433pt, outer sep=0pt,
	fill=gray},
	c3/.style = {trapezium, 
	trapezium angle = 60,
	minimum height= 0.433cm, 
	minimum width= 1.5cm, 
	trapezium stretches=true, 
	inner xsep=1.0pt, inner ysep=0.433pt, outer sep=0pt,
	fill=gray},
	meter/.style ={
	matrix of nodes,
	nodes in empty cells,
	rounded corners,
	draw = black!70,
	fill = black!20,
	inner sep=1pt,
	nodes = {circle, text width=2.5mm, inner sep=0pt, anchor=center, draw=black},
	column 1/.style={nodes={fill=red}},
	column 2/.style={nodes={fill=green}},
	column sep=0.5mm}
}

\biboptions{authoryear}

\begin{document}

\begin{frontmatter}

\title{Queue replacement principle for corridor problems \\ with heterogeneous commuters}

\author[1]{Takara Sakai\corref{cor1}}
\ead{sakai.t.av@m.titech.ac.jp}
\author[2]{Takashi Akamatsu\corref{cor1}}
\ead{akamatsu@plan.civil.tohoku.ac.jp}
\author[3]{Koki Satsukawa\corref{cor1}}
\ead{satsukawa@staff.kanazawa-u.ac.jp}

\address[1]{Department of Civil and Environmental Engineering, Tokyo Institute of Technology, \\2-12-1 W6-9, Ookayama, Meguro, Tokyo 152-8550, Japan.}
\address[2]{Graduate School of Information Sciences, Tohoku University, \\ 6-6 Aramaki Aoba, Aoba-ku, Sendai, Miyagi 980-8579, Japan.}
\address[3]{Institute of Transdisciplinary Sciences for Innovation, Kanazawa University, \\ 2B611 Natural Science \& Technology Hall, Kakuma-machi, Kanazawa, Ishikawa 920-1192, Japan.}

\cortext[cor1]{Corresponding author}

\begin{abstract}
  This study investigates the theoretical properties of a departure time choice problem
  considering commuters' heterogeneity with respect to the value of schedule delay in corridor networks. 
  Specifically, we develop an analytical method to solve the dynamic system optimal (DSO) and dynamic user equilibrium (DUE) problems.
  To derive the DSO solution, we first demonstrate the bottleneck-based decomposition property, i.e.,
  the DSO problem can be decomposed into multiple single bottleneck problems.
  Subsequently, we obtain the analytical solution by applying the theory of optimal transport to each decomposed problem and derive optimal congestion prices to achieve the DSO state. 
  To derive the DUE solution, we prove the queue replacement principle (QRP) that the time-varying optimal congestion prices are equal to the queueing delay in the DUE state at every bottleneck.
  This principle enables us to derive a closed-form DUE solution based on the DSO solution.
  Moreover, as an application of the QRP, we prove that the equilibrium solution under various policies (e.g., on-ramp metering, on-ramp pricing, and its partial implementation) can be obtained analytically.
  Finally, we compare these equilibria with the DSO state.
\end{abstract}


\end{frontmatter}



\section{Introduction}\label{sec:Introduction}

\subsection{Background}
The departure time choice problem proposed by~\cite{Vickrey1969-ic} describes time-dependent traffic congestion as the dynamic user equilibrium (DUE) of commuters' departure time choices.
Many studies have analyzed the departure time choice problem in the classical framework with a single bottleneck and homogeneous commuters, and they have clarified theoretical properties, such as existence, uniqueness and the closed-form analytical solution~\citep{Vickrey1969-ic,Hendrickson1981-cu,Smith1984-wr,Daganzo1985-ls,Newell1987-wk,Arnott1993-mq,Lindsey2004-aw}. 
These results provide important insights for designing traffic management policies, such as dynamic pricing and on-ramp metering (see~\citet{Li2020-zz} for a recent comprehensive review). 
\par
An important property of the single bottleneck model, assuming homogeneous commuters, is that when time-varying congestion prices mimicking the DUE queueing delay pattern are implemented, bottleneck congestion is eliminated without changing the arrival time at the destination of each commuter~\citep{Vickrey1969-ic,Arnott1998-of}.
Such a dynamic pricing pattern leads the traffic state to a dynamic system optimal (DSO) state in which the total system cost is minimized~\citep{Hendrickson1981-cu,Arnott1993-mq,Arnott1994-dy,Laih1994-hi,Akamatsu2021-zg}. 
Therefore, the optimal pricing pattern can be obtained by observing the queueing delay pattern in the DUE state.
Conversely, the queuing delay pattern in the DUE state can be obtained from the optimal pricing pattern in the DSO state~\citep{Iryo2007-ne,Akamatsu2021-zg}.
This means that the DUE flow pattern can be derived using the optimal pricing pattern in the DSO state.
\par
\begin{figure}[tbp]
	\center
  \includegraphics[clip, width=1.0\columnwidth]{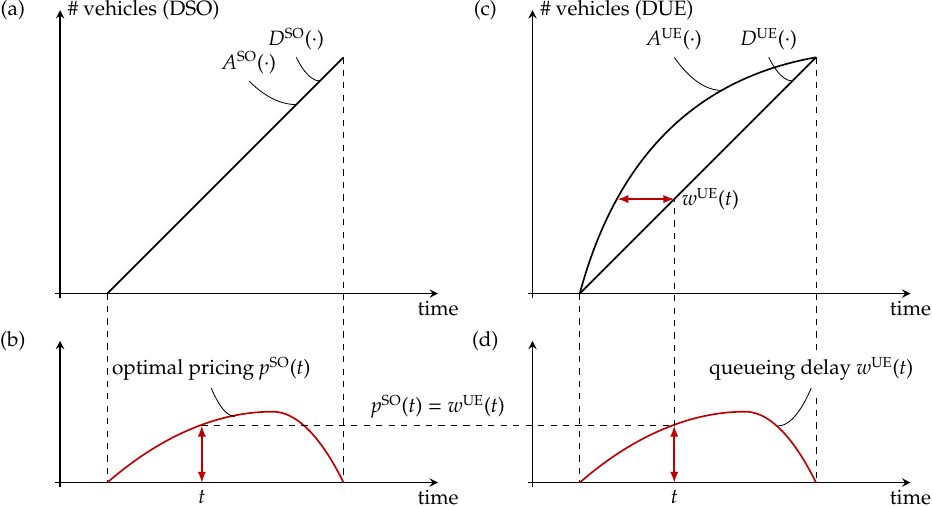}
		\caption{Queue replacement principle in a single bottleneck.}
		\label{fig:SingleQRP}
\end{figure}

\cref{fig:SingleQRP} illustrates this fact by comparing the DSO and DUE states.
  The horizontal axes in \cref{fig:SingleQRP}(a)-(d) represent the arrival time $t \in \ClT$ at the destination of commuters, where $\ClT$ is the morning rush-hour.
\cref{fig:SingleQRP}(a) represents the cumulative arrival curve $\S{A}(\cdot)$ and cumulative departure curve $\S{D}(\cdot)$ for the DSO state. 
Because a queue is eliminated in the DSO state, if the free-flow travel time is ignored, the arrival curve $\S{A}(\cdot)$ is the same as the departure curve $\S{D}(\cdot)$.
\cref{fig:SingleQRP}(b) shows the optimal pricing pattern $\S{p}(t)$, which can be obtained as the optimal Lagrangian multiplier of a linear programming (LP) representing the DSO problem.
\cref{fig:SingleQRP}(d) shows that this optimal pricing pattern equals the queueing delay pattern $\E{w}(t)$ in the DUE state.
In addtion, \cref{fig:SingleQRP}(c) shows that the departure curve in the DUE state $\E{D}(\cdot)$ is the same as that in the DSO state.
Based on $\E{D}(\cdot)$ and $\E{w}(t)$, the arrival curve in the DUE state $\E{A}(\cdot)$ can finally be obtained as shown in \cref{fig:SingleQRP}(c).
\par
	We refer to this remarkable replaceability between the queueing delay and optimal pricing as the queue replacement principle (QRP).
\begin{dfn}
		When the dynamic pricing pattern $\{ \S{p}(t) \}_{t \in \ClT}$ in the DSO state is equal to the queueing delay pattern $\{ \E{w}(t) \}_{t \in \ClT}$ in the DUE state, the QRP holds.
\end{dfn}
\noindent
	The QRP does not only contribute to obtaining the analytical solution, but it also clarifies the efficiency and equity of an optimal pricing scheme from the perspective of welfare analysis.
	In particula, the QRP shows that the travel costs of all commuters do not change with/without implementing optimal dynamic pricing. 
	This means that the QRP is useful for designing an efficient and equitable traffic management scheme.
\par
Considering these contributions of the QRP, it plays a pivotal role in analyzing the theoretical properties in more general settings, such as multiple-bottleneck networks and heterogeneous commuters.
For this problem, \citet{Fu2022-nl} investigated whether QRP holds for the departure time choice problem in corridor networks with multiple bottlenecks, and they demonstrated that QRP holds under certain assumptions. 
However, little is known about whether the QRP holds in the presence of heterogeneous commuters.
  
\subsection{Purpose}
This study proves that the QRP holds for corridor problems with heterogeneous commuters who have different values of schedule delay.
  We prove this QRP condition using the following two-step approach: we derive the analytical (I) DSO and (II) DUE solutions. 
  In part (I), we first formulate the DSO problem as an LP.
  We demonstrate that the DSO solution is established according to the bottleneck-based decomposition property. 
  This property enables us to decompose the DSO problem with multiple bottlenecks into independent single bottleneck problems, which are analytically solvable with the theory of optimal transport
   ~\citep{Rachev1998-bb}. 
   From this analytical solution, we can derive the optimal congestion pricing pattern that achieves the DSO state.
\par
  In part (II), we investigate whether the queueing delay pattern is equivalent to the optimal pricing pattern.
  First, we formulate the DUE problem as a linear complementarity problem (LCP).
  Subsequently, we verify that the queueing delay pattern satisfies the physical requirements of a queue and that the associated DUE flow pattern can be constructed using this queuing delay pattern under a certain assumption which is related to the schedule delay cost function and nonnegative flow condition.
  From this approach, we find that there exists a DUE solution whose queueing delay pattern equals the optimal pricing pattern, i.e., we complete the proof of the QRP under the assumption.
\par
  This QRP implies a replaceability between commuters' pricing and queueing costs at all bottlenecks.
We clarify that this replaceability can independently hold at each bottleneck.
Specifically, we focus on the case in which congestion pricing is only introduced for some bottlenecks.
We demonstrate that the associated equilibrium can be derived by suitably replacing the queueing and pricing costs at each bottleneck. 
Furthermore, as an application of the obtained results, we show that such replaceability of commuters' cost holds when on-ramp-based policies are implemented.
Specifically, we consider equilibriums under on-ramp metering and pricing. 
We reveal that, at each on-ramp, the optimal pricing pattern equals the queueing pattern created by metering in equilibrium, just like the QRP.
Finally, we investigate efficient policies by comparing equilibrium under the bottleneck-based and on-ramp-based policies.

\subsection{Literature review}
For multiple-bottleneck problems, ~\citet{Kuwahara1990-nk} first analyzed the DUE problem in a two-tandem bottleneck network. 
\citet{Kuwahara1990-nk} showed a spatio-temporal sorting property of commuters' arrival times.
In Y-shaped networks with tandem bottlenecks, ~\citet{Arnott1993-mq} theoretically demonstrated a capacity-increasing paradox, and ~\citet{Daniel2009-la} demonstrated a similar paradox in a laboratory setting.
~\citet{Lago2007-wq} studied a similar problem with spillover and merging effects.
  However, applying their approach to analyze cases where an arbitrary number of bottlenecks exist is challenging because they employed the proof-by-cases method, as reported by~\citet{Arnott2011-rb}.
  To overcome this limitation,~\citet{Akamatsu2015-ip} proposed a transparent formulation of the DUE problem in corridor networks with multiple bottlenecks.
  They introduced arrival-time-based variables (i.e., Lagrangian coordinate approach) that facilitate the analysis in corridor networks.\footnote{
  ~\citet{Arnott2011-rb,DePalma2012-la,Li2017-xg}	examined continuum corridor problems.
  Unlike the bottleneck congestion, they focused on "flow congestion" and showed the relationship between velocity and density using the LWR-like traffic flow model.}
  Using the Lagrangian coordinate approach,~\citet{Osawa2018-hg} derived the solution to the DSO problem as part of the long-term location choice problem.
\footnote{
     ~\citet{Osawa2018-hg} considered commuters' heterogeneity; however, they only discussed the first-best traffic flow pattern (DSO problem) because their primary purpose was to obtain long-term policy implications.}
 ~\citet{Fu2022-nl} successfully derived the analytical solution to DUE problems in corridor networks with homogenous commuters. 
\par
	For commuters' heterogeneity, many studies focused on the single bottleneck problem. 
These studies are classified into two categories, depending on how heterogeneity is considered.\footnote{
      There are a few exceptions. 
      ~\citet{Newell1987-wk,Lindsey2004-aw,Hall2018-lg,Hall2021-xo} belong to both categories because they simultaneously consider two types of commuters' heterogeneity.
}
One considered the heterogeneity of the preferred arrival time at the destination as analyzed in~\citet{Hendrickson1981-cu,Smith1984-wr}; and \citet{Daganzo1985-ls}.
They proved the existence and uniqueness of the DUE solution and showed a regularity of the flow pattern, which is called the first-in-first-work principle.
The other category considered the heterogeneity of the value of time;
this includes~\citet{Arnott1988-pl,Arnott1992-jg,Arnott1994-dy,Van_den_Berg2011-qb,Liu2015-ay}; and \citet{Takayama2017-gx}.
They presented welfare analysis and showed that optimal policies can cause inequity in some cases.
\par
  In contrast to these studies, our study investigates the departure time choice problems in the corridor networks with commuters' heterogeneity concerning the value of schedule delay.
  Our contributions are summarized as follows:
\begin{description} 
  \item[(1) Derivation of the analytical DSO solution:]{$ $}
  \par
  This study constructs a systematic approach to solving the DSO problem.
  Specifically, by combining the bottleneck-based decomposition property and the theory of optimal transport, we show that the DSO problem can be solved by sequentially solving single bottleneck problems.
  \item[(2) Proof of the QRP:]{$ $}%
	\par 
  This study investigates the DSO and DUE problems in corridor networks considering commuters' heterogeneity in terms of the value of schedule delay. 
  We successfully prove that the QRP holds under certain conditions.
  \item[(3) Derivation of the analytical DUE solution:]{$ $}%
	\par
	We derive the analytical DUE solution based on the DSO solution and QRP. 
	Moreover, this contribution implies that the solution to an LCP (DUE problem), which is significantly difficult to solve analytically~\citep{Arnott2011-rb,Akamatsu2015-ip}, can be obtained from the solution to the LP (DSO problem). 
  This is a theoretically/mathematically remarkable finding, representing a substantial advancement in the theory of dynamic traffic assignment problems.
  \item[(4) Derivation of the equilibria under various policies using the QRP:]{$ $}%
	\par
  As an application of Contributions (1)-(3), we show that the equilibria under on-ramp metering and on-ramp pricing can be derived using the QRP. This fact clarifies the theoretical relationships between bottleneck-based pricing and on-ramp-based policies. 
  Moreover, we focus on cases in which policies are implemented at certain parts of the network and derive the associated equilibrium solution using the QRP. 
  This means that the QRP enables us to analyze and characterize such an equilibrium state that is neither the pure DUE nor the DSO state.
\end{description}

\subsection{Structure of this paper}
The remainder of this paper is organized as follows:
\cref{sec:Model_settings} introduces the network structure and the heterogeneity of commuters. 
\cref{sec:DSO} presents the formulation of the DSO problem and constructs the systematic approach for deriving its analytical solution. 
In \cref{sec:DUE}, through the proof of the QRP, we derive the analytical solution to the DUE problem using the analytical DSO solution. 
Finally, \cref{sec:Application} presents the application of the QRP. 
\cref{sec:ConcludingRemarks} concludes this paper.
In this paper, all proofs of propositions are given in the appendix.

\begin{figure}[tbp]
  \begin{minipage}[b]{0.62\columnwidth}
			\center
    \includegraphics[clip, width=1.0\columnwidth]{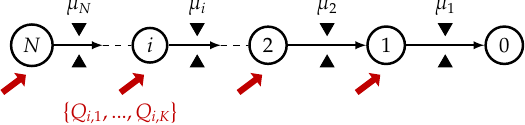}
			\subcaption{Corridor network with $N$ origins and the single destination $0$.}
			\label{fig:CorridorNetwork}
	\end{minipage}
	\begin{minipage}[b]{0.35\columnwidth}
		\center
      \includegraphics[clip, width=0.8\columnwidth]{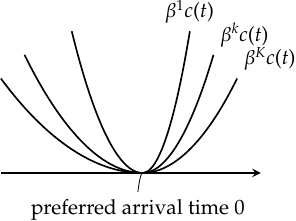}
					\subcaption{Schedule delay cost functions.}
			\label{fig:ScheduleDelayCostFunction}
	\end{minipage}
	\caption{Network and schedule delay cost functions.}
	\label{fig:NetworkAndSDF}
\end{figure}

\section{Model settings}\label{sec:Model_settings}
\subsection{Networks and commuters}
We consider a freeway corridor network consisting of $N$ on-ramps (origin nodes) and a single off-ramp (destination node).
These nodes are numbered sequentially from the destination node, denoted ``$0$'', to the most distant origin node $N$ (\cref{fig:CorridorNetwork}).
The set of origin nodes is denoted by $\ClN \equiv \{ 1,2, \ldots, N \}$.
The link connecting node $i$ to node $(i-1)$ is referred to as link $i$. 
Each link consists of a single bottleneck segment and free-flow segment.
We refer to the bottleneck on link $i$ as bottleneck $i$.
\par
Let $\mu_{i}$ be the capacity of bottleneck $i$.
We also define $\overline{\mu}_{i} \equiv \mu_{i} - \mu_{(i+1)}$, where $\mu_{(N+1)} \equiv 0$.
A queue is formed at each bottleneck when the inflow rate exceeds the bottleneck capacity.
The queue dynamics are modeled by the standard point queue model along with the first-in-first-out (FIFO) principle.
The free-flow travel time from origin (and bottleneck) $i \in \ClN$ to the destination is denoted by $d_{i}$.
\par
From each origin $i$, $Q_{i}$ commuters enter the network and reach their destination during the morning rush-hour $\ClT$.
Commuters are treated as a continuum, and the total mass $Q_{i}$ at each origin $i$ is a given constant.
The commuters are also classified into a finite number $K$ of homogeneous groups with respect to the value of schedule delay. 
The index set of groups is $\ClK = \{1,2,\ldots, K \}$.
The mass of commuters in group $k$ departing from origin $i$ is denoted by $Q_{i, k}$; therefore, $\sum_{k\in\ClK} Q_{i, k} = Q_{i}$.
Commuters in group $k$ departing from origin $i$ are referred to as $(i, k)$-commuters.
\par
The trip cost for each commuter is assumed to be additively separable into free-flow travel, queueing delay, and schedule delay costs.
The schedule delay cost is defined as the difference between the actual and preferred arrival times at the destination.
We assume that all commuters have the same preferred arrival time $t^{\mathrm{p}}=0$.
The schedule delay cost for group $k$ commuters who arrive at the destination at time $t$ is represented by $c^{k}(t)$.
	We assume that $c^{k}(t)$ is represented by $\beta^{k}c(t)$, where $c(t)$ is a base schedule delay cost function with $c(t^{\mathrm{p}})=0$, which is assumed to be strictly quasi-convex and piecewise differentiable $\forall t \in \ClT$, as shown in \cref{fig:ScheduleDelayCostFunction}.
  In addtion, we assume $\mathrm{d} c(t) /\mathrm{d} t > -1$, $\forall t \in \ClT$ as in \citet{Daganzo1985-ls,Lindsey2004-aw}.
$\beta^{k} \in (0, 1]$ is a parameter representing 
the heterogeneity of commuters with $\beta^{K}<... <\beta^{1}=1$.
We also let $\overline{\beta}^{k} \equiv \beta^{k} - \beta^{(k+1)}$, where $\beta^{(K+1)} \equiv 0$.
\par
\begin{figure}[tbp]
	\center
    \includegraphics[clip, width=0.7\columnwidth]{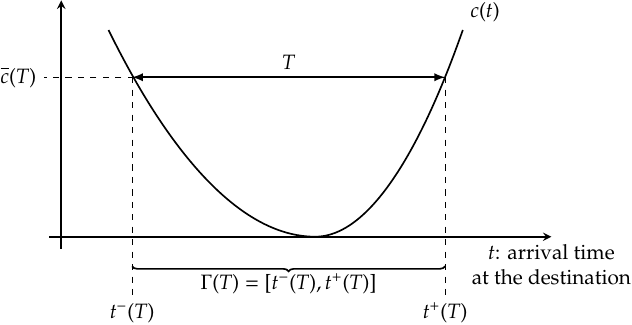}
		\caption{$T$ and $\Gamma(T)$, in relation to $c(t)$.}
    \label{fig:InverseSDF}
\end{figure}

Because $c(t)$ is strictly convex, for a given time length $T$, we can define $t^{-}(T)$ as the unique solution to the equation $c(t) = c(t+T)$. We define $t^{+}(T)$ as $t^{+}(T) = t^{-}(T) + T$.
The base schedule delay cost at $t^{-}(T)$ (and $t^{+}(T)$) is denoted by $\overline{c}(T)$ (i.e., $\overline{c}(T)=c(t^{-}(T))=c(t^{+}(T))$).
We also define time set $\Gamma (T)$ as
\begin{align}
  \Gamma(T)\equiv[t^{-}(T), t^{+}(T)],
  \label{eq:Gamma}
\end{align}
where $\Gamma(T)$ is uniquely determined for $T$ (see Figure~\ref{fig:InverseSDF}).

\subsection{Formulation in a Lagrangian-like coordinate system}
  We primarily describe traffic variables in a \textit{Lagrangian-like coordinate system}~\citep{Kuwahara1993-vt,Akamatsu2015-ip,Akamatsu2021-zg}.
  In this system, variables are expressed in association with the arrival time at the destination, and not at the origin or bottleneck. 
  Such an expression is suitable for considering the ex-post travel time of each commuter during their trip. 
  Therefore, we can easily trace the time-space paths of commuters.
\par
Several variables are introduced for the expression of within-day dynamic traffic flow in this system.
We define $\tau_{i}(t)$ and $\sigma_{i}(t)$ as the arrival and departure times at bottleneck $i\in\mathcal{N}$, respectively, for commuters whose destination arrival time is $t$.
The queuing delay at bottleneck $i$ for the commuters arriving at time $t$ at the destination is denoted by $w_{i}(t)$.
 These variables satisfy the following relationships:
\begin{align}
  &\tau_{i}(t) = t - \sum_{j; j \leq i} w_{j}(t) - d_{i}
	&&\forall i \in \ClN, \quad \forall t \in \ClT,
	\label{eq:tau=s-sumw}
	\\
	&\sigma_{i}(t) = t - \sum_{j; j < i} w_{j}(t) - d_{(i-1)}
	&&\forall i \in \ClN, \quad \forall t \in \ClT.
	\label{eq:sigma=s-sumw}
\end{align}
Note that $\tau_{i}(t)$ must satisfy the following relationship:
\begin{align}
    &\cfrac{\mathrm{d}\tau_{i}(t)}{\mathrm{d}t}\equiv \dot{\tau}_{i}(t) > 0,
    &&\forall i \in \ClN, \quad \forall t \in \ClT,
    \label{eq:TauMonotone}
\end{align}
  where an overdot denotes the derivative of the variable with respect to the destination arrival time $t$.
  This condition guarantees the Lipschitz continuity of cumulative arrival flows, which means that the traffic flow is physically consistent~\citep{Akamatsu2015-ip,Fu2022-nl,Sakai2022-vm}.
Hereafter, we refer to Condition \eqref{eq:TauMonotone} as the consistency condition.

\subsection{Physical conditions for a dynamic traffic flow}
  We formulate physical conditions for dynamic traffic flow.
First, the inflow from each origin must satisfy the following demand conservation condition:
\begin{align}
    &\int_{\ClT} q_{i, k}(t) \mathrm{d} t = Q_{i, k}
	&&\forall k \in \ClK, \quad \forall i \in \ClN,
	\label{eq:DemandCnsvCondition}
\end{align}
where $q_{i,k}(t)\geq 0$ is the inflow rate to the network of $(i, k)$-commuters whose destination arrival time is $t$.
\par
Subsequently, we describe the queueing congestion at each bottleneck.
Let $x_{i}(t)$ be the departure flow rate at link $i\in\mathcal{N}$.
Then, the queueing condition is described as the following complementarity condition:
\begin{align}
    &\begin{dcases}
    x_{i}(\sigma_{i}(t))
    =
    \mu_{i}
    \quad &\mathrm{if}\quad
    w_{i}(t) > 0
    \\
    x_{i}(\sigma_{i}(t))
    \leq
    \mu_{i}
    \quad &\mathrm{if}\quad
    w_{i}(t) = 0
  \end{dcases}
	&&
  \quad \forall k \in \ClK,
  \quad \forall i \in \ClN,
  \quad \forall t \in \ClT.
	\label{eq:queueing_x}
\end{align}

Thirdly, we formulate the flow conservation condition at each node.
This condition requires that for each time, the inflow rate to a node should equal the outflow rate from the node.
Mathematically, this is described as
\begin{align}
    &\cfrac{\mathrm{d}D_{(i+1)}(\sigma_{(i+1)}(t))}{\mathrm{d}t} + \sum_{k\in\mathcal{K}}q_{i,k}(t) = \cfrac{\mathrm{d}A_{i}(\tau_{i}(t))}{\mathrm{d}t}
    &&
    \quad \forall i\in\mathcal{N},
    \quad \forall t\in\mathcal{T}.
\end{align}
Here, the FIFO principle is expressed as
\begin{align}
    &A_{i}(\tau_{i}(t)) = D_{i}(\sigma_{i}(t))
    &&
    \quad \forall i\in\mathcal{N},
    \quad \forall t\in\mathcal{T}.
\end{align}
  Substituting this into the flow conservation condition, we obtain
\begin{align}
    &x_{(i+1)}(\sigma_{(i+1)}(t))\dot{\sigma}_{(i+1)}(t) + \sum_{k\in\mathcal{K}}q_{i,k}(t) =
    x_{i}(\sigma_{i}(t))\dot{\sigma}_{i}(t)
    &&
    \quad \forall i\in\mathcal{N},
    \quad \forall t\in\mathcal{T}.
    \label{eq:flow_conservation_condition_before}
\end{align}
  We have the following equation by recursively applying this toward descendants of each node $i$:
\begin{align}
    \sum_{j;j \geq i}\sum_{k\in\mathcal{K}}q_{i,k}(t) = x_{i}(\sigma_{i}(t))\dot{\sigma}_{i}(t)
    &&
    \quad \forall i\in\mathcal{N},
    \quad \forall t\in\mathcal{T}.
    \label{eq:flow_conservation_condition_after}
\end{align}

Substituting the flow conservation condition \eqref{eq:flow_conservation_condition_after} into \eqref{eq:queueing_x},
we can describe the queueing condition by using the variable $q_{i,k}(t)$, as follows:
\begin{align}
    &\begin{dcases}
    \sum_{j;j \geq i} \sum_{k \in \ClK} q_{j, k}(t)
    = \mu_{i} \dot{\sigma}_{i}(t)
    \quad &\mathrm{if}\quad
    w_{i}(t) > 0
    \\
		\sum_{j;j \geq i} \sum_{k \in \ClK} q_{j, k}(t)
    \leq \mu_{i} \dot{\sigma}_{i}(t)
    \quad &\mathrm{if}\quad
    w_{i}(t) = 0
   \end{dcases}
	 &&\forall i \in \ClN,
   \quad \forall t \in \ClT.
   \label{eq:QueueCondRevised}
\end{align}

\section{Dynamic system optimal problem}\label{sec:DSO}
In this section, an analytical solution to the DSO problem is derived.
\cref{subsec:DSO_Formulation} formulates the DSO problem as an LP.
\cref{subsec:DSO_Aggregate} derives the aggregated DSO flow by decomposing the problem into sub-problems with respect to bottlenecks under a certain assumption. 
  \cref{subsec:DSO_Disaggregate} presents the approach for deriving the disaggregated DSO flow from the aggregate arrival DSO flow by applying the theory of optimal transport~\citep{Rachev1998-bb}.

\subsection{Formulation of the DSO problem} \label{subsec:DSO_Formulation}
We define a DSO state as a state in which the total transport cost is minimized without queues, i.e., congestion externalities are completely eliminated.
Mathematically, the DSO problem is formulated as the following LP~\citep[][]{Osawa2018-hg,Fu2022-nl}: 
\begin{align}	
  \problemname{DSO}
  \\
    \min_{q_{i, k}(t)  \geq \Vt0}.
    \quad
    &
		\sum_{i \in \mathcal{N}} \sum_{k \in \ClK}
    \int_{\ClT} 
    \left( \beta^{k}c(t) + d_{i} \right) 
    q_{i, k}(t) \mathrm{d} t
    &&
    \\
    \mbox{s.t.} \quad
    &  \int_{\ClT} q_{i, k}(t) \mathrm{d} t  = Q_{i, k}
    &&\forall k \in \ClK,
    \quad
      \forall i \in \mathcal{N},
			\label{eq:DSO_ODcncv}
    \\
    & \sum_{j ; j \geq i} q_{j}(t) \leq \mu_{i}
    &&\forall i \in \mathcal{N},
    \quad \forall t \in \ClT,
		\label{eq:DSO_LinkFlow_BNCapa}
\end{align}
where $q_{i}(t)$ is the aggregated destination arrival flow at $t$ of commuters departing from origin $i$:  $q_{i}(t) \equiv \sum_{k \in \ClK} q_{i, k}(t)$.
The objective function is the total schedule delay costs of all commuters.
The first constraint is the flow conservation condition equivalent to \eqref{eq:DemandCnsvCondition}.
The second constraint is the queueing condition equivalent to \eqref{eq:QueueCondRevised} when there are no queues.
\par
The optimality conditions of [DSO] are given as follows~\citep{Luenberger1997-la,Akamatsu2021-zg}: 
    \begin{align}
		&
		0 \leq \S{q}_{i, k}(t) \perp 
		\bracket{ \beta^{k} c(t) + \sum_{j; j\leq i} \S{p}_{j}(t) + d_{i} - \S{\rho}_{i, k}}
		\geq 0
		&&\forall k \in \ClK,
		\quad \forall i \in \mathcal{N},
		\quad \forall t \in \ClT,
		\label{eq:DSO-OC-DTC}
		\\
		&
		0 \leq \S{p}_{i}(t) \perp 
		\bracket{ \mu_{i} - \sum_{j;j \geq i } \S{q}_{j}(t) }
		 \geq 0
		&&\forall i \in \mathcal{N},
		\quad \forall t \in \ClT,
		\\
		& 
		0 \leq \S{\rho}_{i, k} \perp 
		\bracket{ \int_{\ClT} \S{q}_{i, k}(t) \mathrm{d} t - Q_{i, k} }
		\geq 0
		&&\forall k \in \ClK,
		\quad
		\forall i \in \mathcal{N},
	\end{align}
where $\rho_{i, k}$ and $p_{i}(t)$ are the Lagrange multipliers for Constraints~\eqref{eq:DSO_ODcncv} and \eqref{eq:DSO_LinkFlow_BNCapa}, respectively.
The superscript SO indicates that the variables are the optimal solution to [DSO]. 
\par
It is worth noting that we can interpret these optimality conditions as equilibrium conditions under an optimal dynamic congestion pricing scheme~\citep{Arnott1990-ta,Laih1994-hi,Lindsey2012-gg,Chen2015-ku}.
This indicates that $\{ \S{p}_{i}(t) \}_{i \in \ClN}$ can be regarded as an optimal pricing pattern, and $\S{\rho}_{i, k} = \S{\tilde{\rho}}_{i, k} - d_{i}$ can be regarded as the equilibrium commuting cost of $(i,k)$-commuters.
  \cref{eq:DSO-OC-DTC} can then be interpreted as the equilibrium condition for commuters' departure time choice condition under the pricing scheme.
This interpretation from the perspective of equilibrium conditions helps us obtain the analytical solution to [DSO].
\footnote{
		Another interpretation of $\{ \S{\Vtp}(t) \}_{t \in \ClT}$ is the market clearing price pattern under a time-dependent tradable bottleneck permit scheme~\citep{Wada2013-li,Akamatsu2017-bi}. 
		The optimality condition \eqref{eq:DSO_LinkFlow_BNCapa} is interpreted as the demand-supply equilibrium condition for the market of tradable bottleneck permits.
}

\subsection{Aggregated DSO flow pattern}\label{subsec:DSO_Aggregate}
In the following sections, we present a systematic approach to obtain an analytical solution to the DSO problem. 
This section obtains the aggregated DSO flow pattern and shows the decomposition property of [DSO] as a first step of this approach.
\par
  As \citet{Fu2022-nl} noted, the DSO flow pattern in the corridor problem may have arbitrariness.
  To exclude these cases and clarify the analysis approach, we introduce the following assumptions regarding flow and pricing patterns in the DSO state. 
\begin{asm}
	For all bottleneck $i \in \ClN$, the following relationship is satisfied:
	\begin{align}
		&\S{p}_{i}(t) > 0
		\quad \Leftrightarrow \quad 
		\S{q}_{i}(t) > 0
		&&\forall i \in \ClN, \quad \forall t \in \ClT.
		\label{eq:Asm_p>0q>0}
	\end{align}
  The time window $\S{\ClT}_{i} \equiv \mathrm{supp} \left( \S{p}_{i}(t) \right) = \mathrm{supp} \left( \S{q}_{i}(t) \right)$ is convex for all $i \in \ClN$.
	\label{asm:p>0q>0}
\end{asm}
\noindent
  This assumption implies that all commuters experience optimal congestion prices, which are non-zero, at all their passing bottlenecks. This assumption also implies that the network has no false bottlenecks, whose optimal prices are always zero.
  In the corridor network with false bottlenecks, the DSO flow pattern has arbitrariness~\citep{Fu2022-nl}.
  Using this assumption, we can exclude the case where the optimal flow patterns of the closed-form solution have arbitrariness, which causes unnecessary complications for the analysis.
\par
  \cref{asm:p>0q>0} may seem restrictive, but it does not limit the situations that can be analyzed.
  A network that satisfies \cref{asm:p>0q>0} can be constructed from any corridor network with arbitrary capacity patterns by applying the algorithm proposed by \citet{Fu2022-nl}.
  The algorithm detects false bottlenecks, whose optimal prices are always zero\footnote{
    Strictly speaking, the algorithm detects the bottlenecks that do not satisfy the condition \eqref{eq:Asm_p>0q>0} in the original network.}, and reduces the false bottlenecks and corresponding origins by aggregating travel demands from the upstream and downstream origins of each false bottleneck (see \citealt[]{Fu2022-nl}, Appendix B).
  Note that, in this reduced network, $\mu_{i} - \mu_{(i+1)}>0$, $\forall i \in \ClN \setminus \{ N \}$ is the necessary condition of \cref{asm:p>0q>0}.
  If $\mu_{i} - \mu_{(i+1)} \leq 0$, the bottleneck capacity constraint condition \eqref{eq:DSO_LinkFlow_BNCapa} and the non-false bottleneck condition cannot be achieved simultaneously.
\par
Under this assumption, we have the following lemma: 
\begin{lem}
	Suppose that \textbf{\cref{asm:p>0q>0}} holds.
	The arrival time window of $i$-commuters, who depart from origin $i \in \mathcal{N} \setminus \{ N \}$, is included in that of $(i+1)$-commuters in the DSO state: 
	\begin{align}
		&\S{\ClT}_{i} \subset \S{\ClT}_{(i+1)}
		&&\forall i \in \ClN \setminus \{ N \}.
		\label{eq:Lem_S_i_subset_S_i+1}
	\end{align}
	\label{lem:S_i_subset_S_i+1}
\end{lem}
\noindent 
This lemma shows the \textit{spatial sorting property} of the DSO flow pattern \citep{Fu2022-nl}.
  Specifically, the first commuters to depart from origin $i$ arrive at the destination later than the first commuters to depart from the immediate upstream origin $i+1$; the last commuters to depart from origin $i$ arrive earlier than the last commuters to depart from origin $i+1$; and the arrival time windows have nested structures.
\par
\begin{figure}[tbp]
	\center
    \includegraphics[clip, width=0.8\columnwidth]{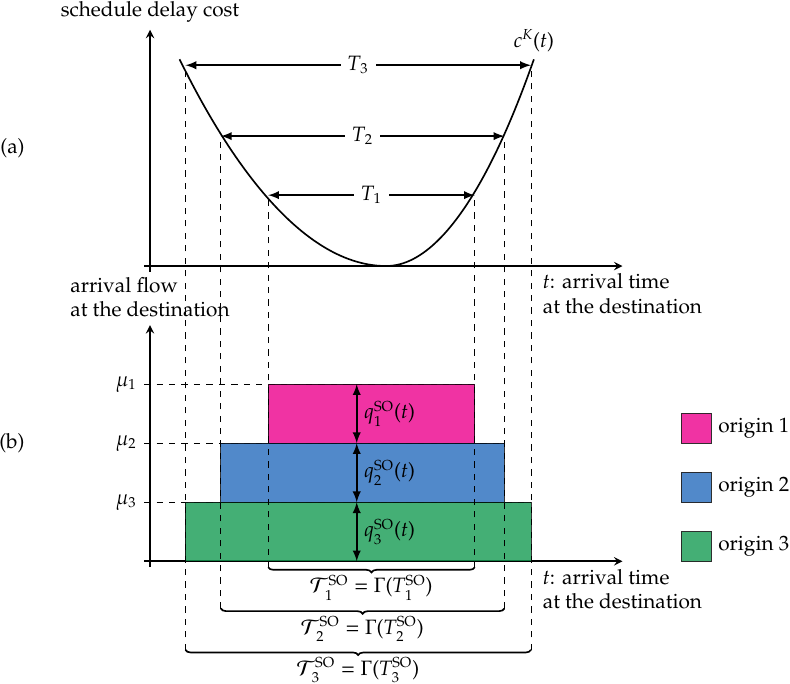}
		\caption{Aggregated DSO arrival flow pattern at the destination with the spatial-sorting property.}
		\label{fig:Spatial_Sorting}
\end{figure}

Using this property, we obtain the aggregated DSO flow pattern $\{ \S{q}_{i}(t) \}$, as follows: 
\begin{lem}
	Suppose that \textbf{\cref{asm:p>0q>0}} holds, 
	the aggregated DSO flow is given as:
	\begin{align}
		\S{q}_{i}(t) 
		&=
		\begin{dcases}
			\overline{\mu}_{i}
			\quad &\mathrm{if}\quad t \in \S{\ClT}_{i}
			\\
			0 \quad &\mathrm{otherwise}
		\end{dcases}
	 &&\forall i \in \ClN,
	 \label{eq:Lem_q=mu}
	\end{align}
    where $\overline{\mu}_{i} = \mu_{i} - \mu_{(i+1)}$ and $\mu_{(N+1)} = 0$.
    The length of the arrival time window $\S{T}_{i} \equiv | \S{\mathcal{T}}_{i} |$ is determined as 
     $\S{T}_{i} = \sum_{k \in \ClK} Q_{i,k} /\overline{\mu}_{i}$.
	\label{lem:q=mu}
\end{lem}
\noindent This lemma shows that the aggregated DSO flow pattern becomes an \textit{all-or-nothing pattern}.
The aggregated inflow rate $\S{q}_{i}(t)$ equals the difference in the capacities between the immediate downstream and upstream bottlenecks during $\S{\mathcal{T}}_{i}$; otherwise, the inflow rate becomes zero.
  \cref{fig:Spatial_Sorting} (a) shows the aggregated DSO flow pattern with $N=3$.
  We find that the spatial sorting property makes the aggregated DSO flow pattern as a staircase.
  \cref{fig:Spatial_Sorting} also illustrates that the time window $\S{\mathcal{T}}_{i}$ is determined by $\Gamma(\S{T}_{i})$, where $\Gamma(\cdot)$ is defined by \cref{eq:Gamma}.
  This fact is derived from \cref{lem:SubSol}, which will be discussed later. 
  From \cref{lem:SubSol}, the first and last commuters from origin $i$ to arrive at the destination are both in group $K$.
  Based on the optimality condition \eqref{eq:DSO-OC-DTC} and \cref{asm:p>0q>0}, 
  we find $\beta^{K}c(\min. \{\S{\mathcal{T}}_{i} \})=\beta^{K}c(\max. \{\S{\mathcal{T}}_{i} \})$ for all $i \in \ClN$.
  Hence, the time window $\S{\mathcal{T}}_{i}$ is determined by $\Gamma(\S{T}_{i})$ as shown in \cref{fig:Spatial_Sorting}(a) and (b).

\begin{figure}[tbp]
	\center
    \includegraphics[clip, width=0.9\columnwidth]{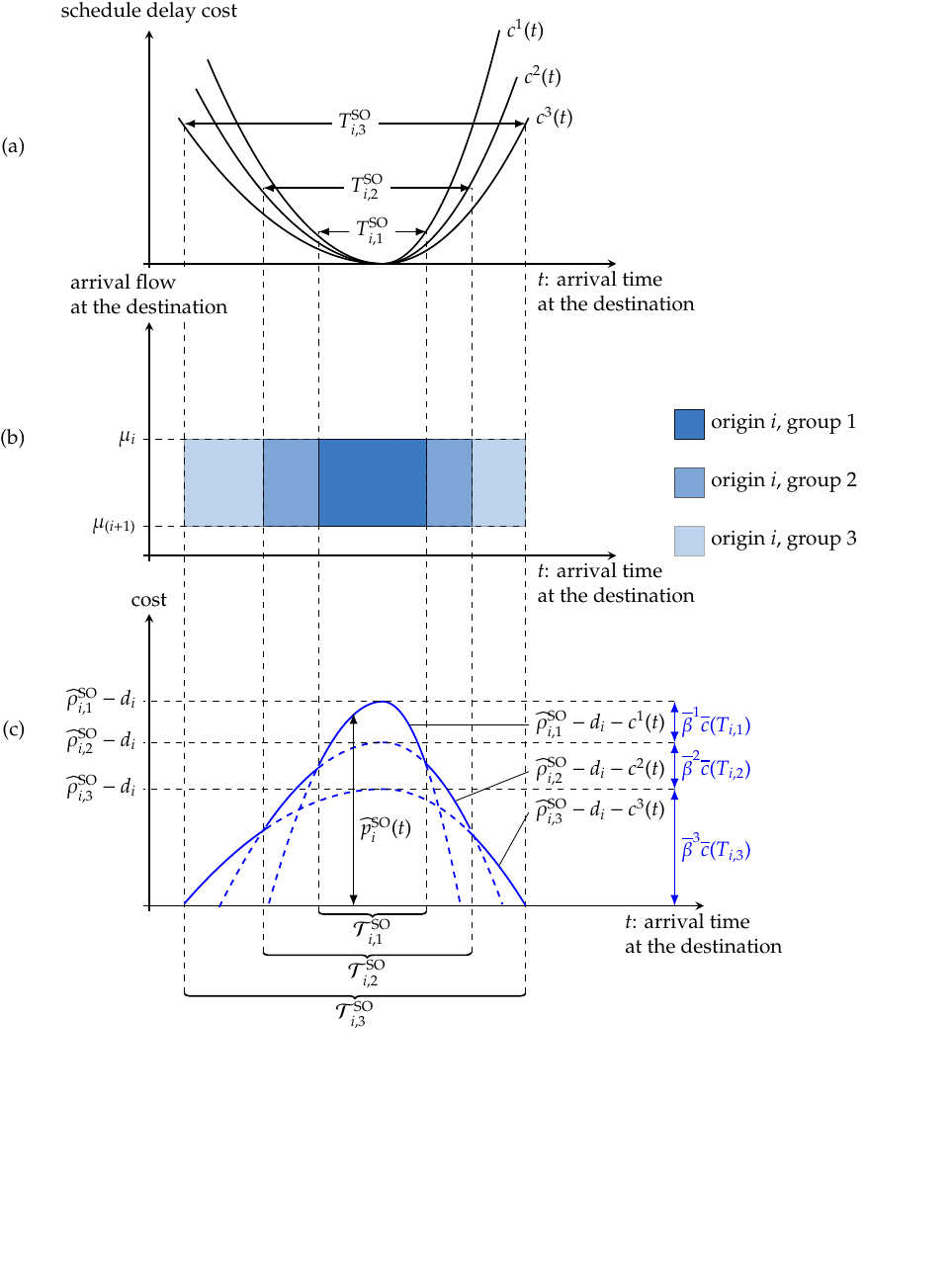}
		\caption{Disaggregated DSO arrival flow pattern at the destination with the temporal-sorting property.}
		\label{fig:Temporal_Sorting}
\end{figure}

Interestingly, \textbf{\cref{lem:q=mu}} is more than a simple contribution toward clarifying the aggregated DSO flow pattern.
This enables us to \textit{decompose the DSO problem} in a corridor network into bottleneck-based sub-problems, which are useful for obtaining the disaggregated DSO flow pattern, i.e., inflow rate of commuters in each group.
The lemma eliminates the interaction of bottlenecks in Constraint \eqref{eq:DSO_LinkFlow_BNCapa}, and Constraint \eqref{eq:DSO_LinkFlow_BNCapa} can be independently rewritten as follows:
\begin{align}
	&\sum_{k \in \ClK} q_{i, k}(t) =  \overline{\mu}_{i}
	&&\forall i \in \ClN, \quad \forall t \in \S{\ClT}_{i}.
	\label{eq:DSO_Disaggregated_DSOFlow_BNCapa}
\end{align}
	Because the RHS of \cref{eq:DSO_Disaggregated_DSOFlow_BNCapa} is constant, this form implies the separation of \cref{eq:DSO_LinkFlow_BNCapa} by each bottleneck $i$.
	Additionally, the objective function of [DSO] can also be separated by each origin (bottleneck) $i$, and the demand conservation constraint \eqref{eq:DSO_ODcncv} is already independent at each origin $i$.
	Consequently, we find that the objective function and all constraints of [DSO] can be separated by each origin (bottleneck) $i$.
	This means we can independently consider each bottleneck-based sub-problem without considering the other sub-problems. 
	Hence, for obtaining the disaggregated DSO flow pattern, we do not have to deal with complex spatial interactions between multiple bottlenecks; we solve the sub-problem, which has the same mathematical structure as a single bottleneck problem. 
	This is summarized in the following lemma:

\begin{lem}[Decomposition Property]
    Suppose that \textbf{\cref{asm:p>0q>0}} holds.
    Then, the solution to [DSO] can be obtained by solving the following sub-problems [DSO-Sub$(i)$] for every bottleneck $i\in\mathcal{N}$:
    \begin{align}
      \problemname{DSO-Sub$(i)$} 
      \\
      \min_{\{ q_{i, k}(t) \}_{k \in \ClK, t \in \ClT} \geq \Vt0}.
      \quad
      &
			\sum_{k \in \ClK}
      \int_{\S{\ClT}_{i}} 
      \left(
        \beta^{k}c(t) + d_{i}
      \right)
       q_{i, k}(t) \mathrm{d} t
      \label{eq:Objective_SubProblem}
      &&
      \\
      \mbox{s.t.}
      \quad
      &\int_{\S{\ClT}_{i}} q_{i, k}(t) \mathrm{d} t
      = Q_{i, k}
      &&\forall k \in \ClK,
			\label{eq:DSO-sub_ODcncv}
      \\
      &\sum_{k \in \ClK} q_{i, k}(t) = \overline{\mu}_{i}
      &&\forall t \in \S{\ClT}_{i}.
			\label{eq:DSO-sub_BNCapa}
    \end{align}
  \label{lem:SubProb}
\end{lem}	
The optimality conditions of [DSO-Sub$(i)$] are given as follows: 
\begin{align}
		&
		0 \leq \S{q}_{i, k}(t) \perp 
		\bracket{\beta^{k} c(t) + \S{\widehat{p}}_{i}(t) + d_{i} - \S{\widehat{\rho}}_{i, k}}
		\geq 0
		&&\forall k \in \ClK,
		\quad \forall t \in \ClT,
		\label{eq:DSO-Sub-OC-DTC}
		\\
		&\overline{\mu}_{i} - \sum_{k \in \ClK} \S{q}_{i,k}(t)
		 = 0
		&&\forall t \in \ClT,
		\\
		&\int_{\ClT} \S{q}_{i, k}(t) \mathrm{d} t - Q_{i, k} = 0
		&&\forall k \in \ClK,
	\end{align}
where $\widehat{\rho}_{i, k}$ and $\widehat{p}_{i}(t)$ are the Lagrange multipliers for constraints~\eqref{eq:DSO-sub_ODcncv} and \eqref{eq:DSO-sub_BNCapa}, respectively.

\subsection{Disaggregated DSO flow pattern}\label{subsec:DSO_Disaggregate}
	The remainder of our approach is to obtain the disaggregated DSO flow pattern by  analytically solving the sub-problems. 
	We remember that this sub-problem [DSO-Sub(i)] has the same structure as the single bottleneck problem.
  Note that total free-flow travel cost is constant and does not affect the optimal flow pattern.
	As discussed in~\cite{Akamatsu2021-zg}, the mathematical structure of the single bottleneck problem is the same as that of the optimal transport problem~\citep{Rachev1998-bb}, which is well known as Hitchcock's transportation problem in operations research and transportation fields. 
	Moreover,~\cite{Akamatsu2021-zg} demonstrated that the analytical solution to the single bottleneck problem could be derived using the theory of optimal transport when the transport cost function, $c^{k}(t)$ in Eq.~\eqref{eq:Objective_SubProblem}, satisfies the \textit{submodularity property}.
	Therefore, it is sufficient to confirm that the schedule delay cost function satisfies that property.
\par
We first introduce the definition of the submodularity property:
\begin{dfn}
	Function $c: \mathbb{R}^{2} \rightarrow \mathbb{R}$ is submodular if, and only if,
	\begin{align}
    c(x, y) + c(x', y') \leq c(x, y') + c(x', y)
    \quad \forall x \leq x', y\leq y'.
		\label{eq:dfn_submodular_1}
  \end{align}
	If the inequality in \cref{eq:dfn_submodular_1} strictly holds for 
	$x' < x$, and $y < y'$, $c$ is called a strict submodular function.
	In addition, if the inequality \cref{eq:dfn_submodular_1}
	holds in the opposite direction, then function $c$ is supermodular.
  \label{dfn:SubModular}
\end{dfn}
\noindent 
  In our sub-problem [DSO-Sub$(i)$], the schedule delay cost function $\beta^{k}c(t)$ is submodular and supermodular in  
  $t \geq 0$ and $t<0$, respectively:
\begin{lem}
    The schedule delay cost function $\beta^{k}c(t)$ 
    is the submodular function for all $t\geq 0$:
  \begin{align}
    \beta^{k}c(t) + \beta^{k+1}c(t') <
    \beta^{k+1}c(t) +  \beta^{k}c(t')
    &&\forall k \in \ClK, \quad \forall t'>t.
    \label{eq:LateSubmodular}
  \end{align}
  \label{lem:LateSubmodular}
\end{lem}
\begin{lem}
    The schedule delay cost function $\beta^{k}c(t)$ 
    is the supermodular function for all $t<0$:
  \begin{align}
    \beta^{k}c(t) + \beta^{k+1}c(t') >
    \beta^{k+1}c(t) +  \beta^{k}c(t')
    &&\forall k \in \ClK, \quad \forall t'>t.
    \label{eq:EarlySupermodular}
  \end{align}
  \label{lem:EarlySupermodular}
\end{lem}
A detailed discussion was given by~\cite{Akamatsu2021-zg}, but here we remark that the submodularity (supermodularity) of $c^{k}(t)$ prioritizes all groups.
Specifically, in the optimal solution, all $(i, k)$-commuters arrive at the destination closer to the preferred arrival time than the $(i, k+1)$-commuters.
In other words, for all origins $i \in \mathcal{N}$, the arrival time window of the $(i, k+1)$-commuters includes that of the $(i, k)$-commuters.
Based on this \textit{temporal sorting property} and the optimality condition of [DSO-Sub$(i)$], we obtain the analytical solution to [DSO-Sub$(i)$] as follows:
\begin{lem}[~\citet{Akamatsu2021-zg}]
	The following $\{ \S{\Vtq}(t) \}_{t \in \ClT}$, $\S{\widehat{\Vtrho}}$, 
	and $\{ \S{\widehat{\Vtp}}(t) \}_{t \in \ClT}$ are solutions to [DSO-Sub$(i)$]:
  \begin{align}
    &\S{q}_{i, k}(t) =
    \begin{dcases}
      \overline{\mu}_{i}
        &\mathrm{if}\quad
      \forall t \in \S{\ClT}_{i, k} \setminus \S{\ClT}_{i, (k-1)}
      \\
      0
        &\mathrm{otherwise}
    \end{dcases}
    &&\forall k \in \ClK,
    \label{eq:inLem_Sub_Sol}
		\\
		&\S{\widehat{\rho}}_{i, k} =
    \sum_{l;l \geq k}
    \overline{\beta}^{l}
    \overline{c}(\S{T}_{i, k}) + d_{i}
    &&\forall k \in \ClK,
    \label{eq:DSO-sub_sol_rho}
    \\
    &\S{\widehat{p}}_{i}(t)
    =
    \S{\widehat{\rho}}_{i, k} - \beta^{k} c(t)
    &&\forall t \in \S{\ClT}_{i, k} \setminus \S{\ClT}_{i, (k-1)},
    \label{eq:DSO-sub_sol_p}
  \end{align}
	  where $\S{\ClT}_{i, k} = \Gamma (\S{T}_{i, k})$
		and $\S{T}_{i, k} = \sum_{l;l\leq k} Q_{i, l} / \overline{\mu}_{i}$.
    Note that $\overline{\beta}^{k} = \beta^{k} - \beta^{(k+1)}$ and $\beta^{(K+1)} = 0$.
  \label{lem:SubSol}
\end{lem}
\par
\cref{fig:Temporal_Sorting} illustrates the solution to the sub-problem [DSO-Sub$(i)$] with $K=3$.
\cref{fig:Temporal_Sorting} (b) represents the disaggregated DSO flow pattern $\{ \S{\Vtq}(t) \}_{t \in \ClT}$ with the temporal sorting property that the arrival time window of the $(i, k+1)$-commuters includes that of the $(i, k)$-commuters.
The arrival time windows are determined by the schedule delay cost functions as shown in \cref{fig:Temporal_Sorting} (a) and (b).
Moreover, \cref{fig:Temporal_Sorting} (c) illustrates the optimal Lagrangian multipliers.
\par
	Combining these analytical solutions to the sub-problems, we immediately construct the complete DSO flow pattern.
	Moreover, the analytical solutions to the sub-problems enable us to derive the optimal Lagrangian multipliers of [DSO-LP], which represent the congestion price and equilibrium commuting cost.
	Specifically, based on the optimality condition of [DSO-LP], we find that $\S{\rho}_{i,k}= \S{\widehat{\rho}}_{i,k}$ and $\S{p}_{i}(t) = \min. \{0, \  \S{\widehat{p}}_{i}(t) - \S{\widehat{p}}_{(i-1)}(t) \}$.
	This is summarized as follows:
\begin{pro}[Solution to the DSO Problem]
	If \cref{asm:p>0q>0} holds, then,
	the following $\{ \S{\Vtq}(t) \}_{t \in \ClT}$ is a solution to [DSO], and 
	the following $\S{\Vtrho}$ and $\{ \S{\Vtp}(t) \}_{t \in \ClT}$
	are optimal Lagrangian multipliers of [DSO]:
\begin{align}
	\S{q}_{i, k}(t) &=
	\begin{dcases}
		\overline{\mu}_{i}
		&\mathrm{if}\  t \in \S{\ClT}_{i, k}
		\setminus \S{\ClT}_{i, (k-1)}
		\\
		0
		&\mathrm{otherwise}
	\end{dcases}
	&&\forall k \in \ClK, \quad \forall i \in \mathcal{N},
	\label{eq:DSO_sol_q}
	\\
	\S{\rho}_{i, k} &=
	\sum_{l;l \geq k}
	\overline{\beta}^{l}
	\overline{c}(\S{T}_{i, k}) + d_{i}
	&&\forall k \in \ClK, \quad \forall i \in \mathcal{N},
	\label{eq:DSO_sol_rho}
	\\
	\S{p}_{i}(t)
	&=
	\begin{dcases}
		\S{\rho}_{i, k} - \beta^{k} c(t) - \sum_{j;j<i} \S{p}_{j}(t)
		&\mathrm{if} \
		t \in \S{\ClT}_{i, k} \setminus \S{\ClT}_{i, (k-1)}
		\\
		0 &\mathrm{otherwise}
	\end{dcases}
	&&\forall i \in \mathcal{N}.
	\label{eq:DSO_sol_p}
\end{align}
\label{pro:DSO_Solution}
\end{pro}
\par
\begin{figure}
	\center
    \includegraphics[clip, width=0.95\columnwidth]{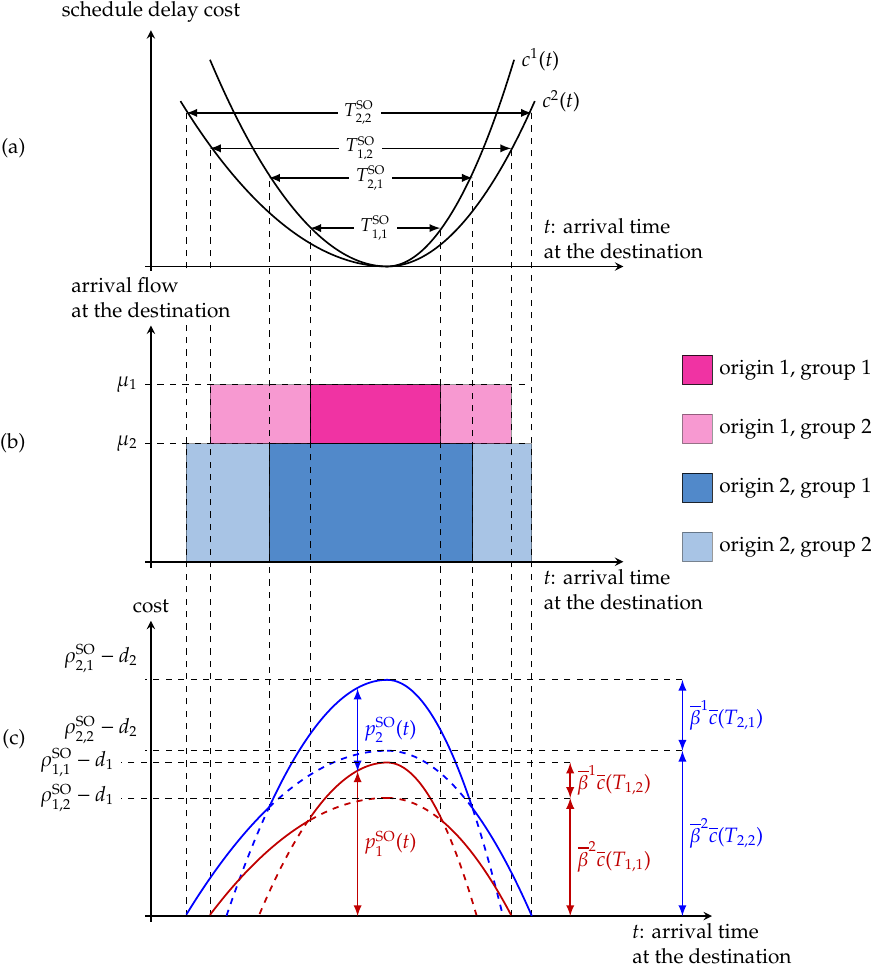}
		\caption{DSO solution.}
		\label{fig:DSO-Solution}
\end{figure}

	\cref{fig:DSO-Solution} illustrates the solution to the DSO problem for the case of $(N,K) = (2,2)$.
  \cref{fig:DSO-Solution} (a) depicts the relationship between the schedule delay cost function and the arrival time window of the commuters.
  \cref{fig:DSO-Solution} (b) and (c) illustrates the DSO arrival flow pattern $\{\S{\Vtq}\}_{t \in \ClT}$ and the optimal Lagrangian multipliers $\{ \S{\Vtp}(t) \}_{t \in \ClT}$, $\S{\Vtrho}$, respectively.
  We find that $\{\S{\Vtq}\}_{t \in \ClT}$, $\{ \S{\Vtp}(t) \}_{t \in \ClT}$, and $\S{\Vtrho}$ can be constructed by combining the solutions to the sub-problems.

\section{Queue replacement principle and the dynamic user equilibrium problem}\label{sec:DUE}
  In this section, we prove the QRP and present an approach to construct a DUE solution using the DSO solution obtained in \cref{sec:DSO}.
\cref{subsec:Formulation_DUE} formulates the DUE problem according to~\cite{Akamatsu2015-ip}.
In \cref{subsec:Proof_QRP}, we prove the QRP and derive the analytical DUE solution.
Finally, \cref{subsec:DSOvsDUE} compares the solutions and demonstrates the essential similarities and differences between the DSO and DUE states.

\subsection{Formulation of the DUE problem}\label{subsec:Formulation_DUE}
  To formulate the DUE problem, we first introduce the commuter's trip cost.
  As mentioned above, the commuter's trip cost is measured based on the arrival time at the destination. 
  It is assumed to be additively separable into schedule delay costs, queuing delay costs, and free-flow travel times.
  The trip cost of $(i,k)$-commuters arriving at time $t$ at the destination is defined as follows:
  \begin{align}
    &C_{i,k}(t) = \beta^{k} c(t) + \alpha 
    \left( \sum_{j;j \leq i} w_{j}(t) + d_{i} \right)
    &&\forall i \in \ClN, \quad \forall k \in \ClK, \quad \forall t \in \ClT,
    \label{eq:TripCost}
  \end{align}
  where $\alpha$ is a parameter representing the value of time. This study assumes $\alpha = 1$ regardless of the commuter's
  group. 
In the DUE state, all commuters can not reduce their commuting costs by unilaterally changing their arrival time.
This is expressed as the following linear complementarity condition:
\begin{align}
	&\begin{dcases}
		\sum_{j; j\leq i} \E{w}_{j}(t) + \beta^{k}c(t) 
    + d_{i} = \E{\rho}_{i, k}
		\quad &\mathrm{if}\quad \E{q}_{i, k}(t) > 0
		\\
		\sum_{j; j\leq i} \E{w}_{j}(t) + \beta^{k}c(t) 
    + d_{i} \geq \E{\rho}_{i, k}
		\quad &\mathrm{if}\quad \E{q}_{i, k}(t) = 0
	 \end{dcases}
    &&\forall k \in \ClK,
    \quad \forall i \in \mathcal{N},
    \quad \forall t \in \ClT,
    \label{eq:DUE-DTC_pre}
\end{align}
  where $\E{w}_{i}(t)$ is a queueing delay cost for commuters arriving at time $t$ at the destination, and $\E{\rho}_{i, k}$ is an equilibrium commuting cost of the $(i, k)$-commuters.
The superscript $\mathrm{UE}$ indicates that the variables are defined in the DUE problem.
\par
  In addition to the departure time choice condition \eqref{eq:DUE-DTC_pre},
  the queueing condition \eqref{eq:QueueCondRevised}, demand conservation condition \eqref{eq:DemandCnsvCondition}, and consistency condition \eqref{eq:TauMonotone}
  in \cref{sec:Model_settings} must be satisfied.
  Thus, the DUE problem is formulated as the following LCP:
\begin{align}
  \problemname{DUE-LCP}
  \\
   \mbox{Find} \quad
	 &
		 \{\E{\Vtq}(t)\}_{t \in \ClT},
		 \{\E{\Vtw}(t)\}_{t \in \ClT}, 
		 \E{\Vtrho}
	 \notag
   \\
   \mbox{such that} \quad
   &
	 0 \leq \E{q}_{i, k}(t) \perp 
	 \bracket{ \sum_{j; j\leq i} \E{w}_{j}(t) + \beta^{k} c(t) + d_{i} - \E{\rho}_{i, k} }
	 \geq 0
	 &&\forall k \in \ClK,
   \quad \forall i \in \mathcal{N},
   \quad \forall t \in \ClT,
	 \label{eq:DUE-DTC}
   \\
   &
	 0 \leq  \E{w}_{i}(t) \perp 
	 \bracket{ \mu_{i} \E{\dot{\sigma}}_{i}(t) - \sum_{j;j\geq i }  \E{q}_{j}(t) }
	 \geq 0
	 &&\forall i \in \mathcal{N},
   \quad \forall t \in \ClT,
   \label{eq:DUE_queue}
	 \\
	 & 0 \leq \E{\rho}_{i,k} \perp
	 \bracket{ \int_{\ClT}  \E{q}_{i, k}(t) \mathrm{d} t - Q_{i, k} }
	 \geq 0
	 &&\forall k \in \ClK,
	 \quad \forall i \in \mathcal{N},
	 \label{eq:DUE_Demand}
		\\
		& 1- \sum_{j;j \leq i} \E{\dot{w}}_{j}(t) > 0 
		&&\forall i \in \ClN, \quad \forall t \in \ClT.
		\label{eq:DUE_queue_consistency}
\end{align}

\subsection{Queue replacement principle and the DUE solution}\label{subsec:Proof_QRP}
We first formally introduce the QRP concept in the multiple-bottleneck networks:
\begin{dfn}[QRP]
    If there exists a DUE state such that the queueing delay pattern coincides with the optimal pricing pattern:
    \begin{align}
    	&\E{w}_{i}(t) = \S{p}_{i}(t)
    	&&\forall i \in \ClN, \quad \forall t \in \ClT,
    \end{align}
    then, the QRP holds.
	\label{dfn:QRP_Multiple}
\end{dfn}
  We explore whether the QRP holds by deriving the equilibrium commuting cost and flow pattern that satisfy the DUE condition when $\{ \S{\Vtp}(t) \}_{t \in \ClT}$ is substituted for $\{ \E{\Vtw}(t) \}_{t \in \ClT}$.
  Consequently, we found that such an equilibrium cost and flow pattern exist if the schedule delay cost function satisfies certain conditions.
  We derive this conclusion from the following approach:
  First, by comparing the DSO problem's optimality condition \eqref{eq:DSO-OC-DTC} and the DUE problem's departure time choice condition \eqref{eq:DUE-DTC}, we observe that $\E{\rho}_{i,k} = \S{\rho}_{i,k}$ for all $(i,k)$-commuters.
  Second, based on this, we conjecture that arrival time window $\E{\ClT}_{i,k} \equiv \{ t \in \ClT \mid \E{q}_{i,k}(t) > 0 \}$ is the same as the arrival time window in the DSO state $\S{\ClT}_{i,k}$ for all $(i,k)$-commuters.
    From this conjecture and the queueing delay condition \eqref{eq:DUE_queue}, we can derive the flow pattern $\{ \Vtq(t) \}_{t \in \ClT}$ as follows:
  \begin{align}
    q_{i, k}(t) &=
      \begin{dcases}
        \mu_{i} \E{\dot{\sigma}}_{i}(t)
        -
        \mu_{(i+1)} \E{\dot{\sigma}}_{(i+1)}(t)
          &\mathrm{if}\quad t \in
          \E{\ClT}_{i, k} \setminus \E{\ClT}_{i, (k-1)}
        \\
        0
        &\mathrm{otherwise}\quad
      \end{dcases}
      &&\forall k \in \ClK, \quad \forall i \in \ClN.
  \end{align}
  However, this flow pattern can be negative, i.e., inconsistent with the departure time choice condition \eqref{eq:DUE-DTC} of the DUE condition.
  Therefore, to prevent the negative flow, we introduce a condition of the schedule delay cost function described below.
  Finally, under this condition, we prove the conjecture by confirming that this flow pattern satisfies the demand conservation conditions.
  \par
    Because this approach implies that there exist an equilibrium commuting cost and flow pattern that satisfy the DUE condition when $\{ \E{\Vtw}(t) \}_{t \in \ClT} = \{ \S{\Vtp}(t) \}_{t \in \ClT}$, we conclude that the QRP holds under the condition of the schedule delay cost function.
  This is summarized in the following theorem and proposition:
  \begin{thm}[Sufficient condition for the QRP]
    Suppose that the schedule delay cost function satisfies the following conditions:
    \begin{align}
      &\max \left\{ 
        -1, \quad - \dfrac{ \mu_{i} - \mu_{(i+1)} }{\mu_{i} \beta^{\max}  - \mu_{(i+1)} \beta^{\min} }
      \right\} < 
      \dot{c}(t) < \dfrac{\mu_{i}- \mu_{(i+1)}}{\beta^{\max}\mu_{(i+1)}}
      &&\forall i \in \ClN \setminus \{N\}, \quad  \forall t \in \ClT,
      \label{eq:QRP_Corridor}
      \\
      &\text{where} \quad \beta^{\min} \equiv \min_{k \in \ClK} \left\{ \beta^{k}\right\}, \quad \beta^{\max} \equiv \max_{k \in \ClK}  \left\{ \beta^{k} \right\}.
      \notag
    \end{align}
    Then, the QRP holds. We refer to the condition \eqref{eq:QRP_Corridor} as the QRP condition.\footnote{
        The condition \eqref{eq:QRP_Corridor} is a generalization of the condition of~\citet{Fu2022-nl}, which assumes homogeneous commuters.
        If all commuters are homogenous, i.e., $\beta^{\max} = \beta^{\min} = 1$, the condition \eqref{eq:QRP_Corridor} can be simplified $-1 < \dot{c}(t) < (\mu_{i}- \mu_{(i+1)})/\mu_{(i+1)}$, which corresponds to the conditions (3.20a) and (3.20b) in ~\citet{Fu2022-nl}.
      }
    \label{thm:QRP_Corridor}
  \end{thm}
  \begin{pro}[Solution to the DUE Problem]
      Suppose that the QRP condition \eqref{eq:QRP_Corridor} holds.
      Then, the following $\{ \E{\Vtq}(t) \}_{t \in \ClT}$, $\{ \E{\Vtw}(t) \}_{t \in \ClT}$ and $\E{\Vtrho}$ are the solutions to [DUE-LCP]:
    \begin{align}
      \E{\ClT}_{i, k} &= \S{\ClT}_{i, k}
      &&\forall k \in \ClK, \quad \forall i \in \ClN.
      \\
      \E{q}_{i, k}(t) &=
      \begin{dcases}
        \mu_{i} \E{\dot{\sigma}}_{i}(t)
        -
        \mu_{(i+1)} \E{\dot{\sigma}}_{(i+1)}(t)
          &\mathrm{if}\quad t \in
          \E{\ClT}_{i, k} \setminus \E{\ClT}_{i, (k-1)}
        \\
        0
        &\mathrm{otherwise}\quad
      \end{dcases}
      &&\forall k \in \ClK, \quad \forall i \in \ClN,
      \label{eq:DUE_final_sol_q}
      \\
      \E{w}_{i}(t) &= \S{p}_{i}(t)
      &&\forall i \in \mathcal{N}, \quad \forall t \in \ClT,
      \label{eq:DUE_final_sol_w}
      \\
      \E{\rho}_{i, k} &= \S{\rho}_{i, k}
      &&\forall k \in \ClK,
      \quad \forall i \in \mathcal{N},
      \label{eq:DUE_final_sol_rho}
      \\
      \mbox{where} &\quad \E{\sigma}_{i}(t) = t - \sum_{j;j\leq i} \E{w}_{j}(t)
      &&\forall i \in \ClN, \quad \forall t \in \ClT.
    \end{align}
    \label{pro:DUE_Solution}
  \end{pro}
  \noindent
    Considering $\beta^{k} \in (0, 1]$ for all $k \in \ClK$, we can simplify the QRP condition \eqref{eq:QRP_Corridor} to make it independent of $\{ \beta^{k} \}$ as follows:
  \begin{align}
    &
      \dfrac{\mu_{(i+1)}}{\mu_{i}} -1 < \dot{c}(t) < \dfrac{\mu_{i}}{\mu_{(i+1)}} -1
    &&
    \forall i \in \ClN \setminus \{ N \},
    \quad \forall t \in \ClT.
    \label{eq:QRP_Corridor_Simplified}
  \end{align}
    This condition \eqref{eq:QRP_Corridor_Simplified} is a sufficient condition for the condition \eqref{eq:QRP_Corridor}.
    Note that the condition \eqref{eq:QRP_Corridor_Simplified} always holds $\dot{c}(t) > -1$ because $\mu_{i+1}/\mu_{i}>0$.

\subsection{Comparison between the DUE and DSO solutions}
\label{subsec:DSOvsDUE}
\begin{figure}[tbp]
	\center
    \includegraphics[clip, width=0.8\columnwidth]{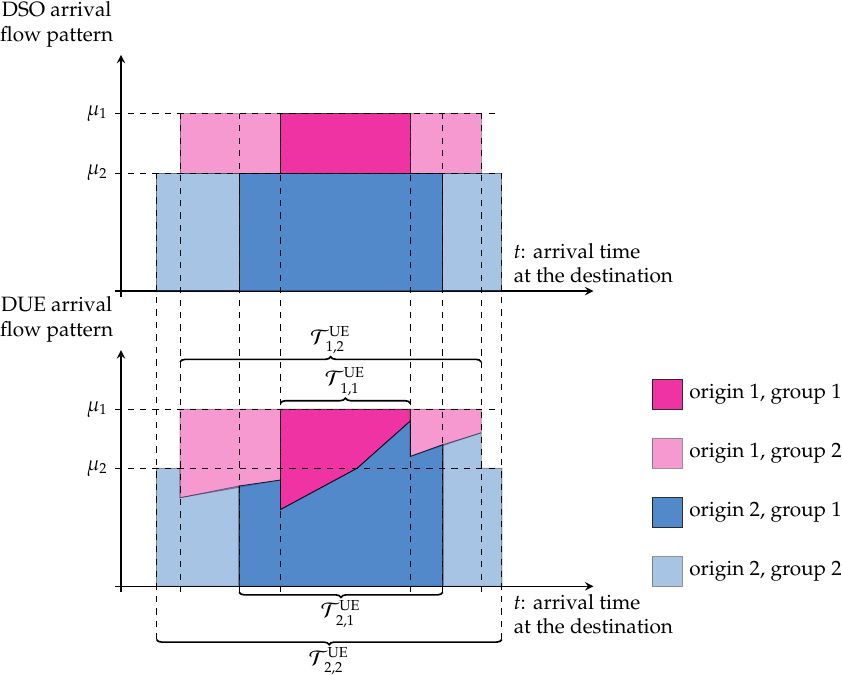}
		\caption{DSO and DUE arrival flow pattern at the destination.}
		\label{fig:DSO-DUE_FlowPtn}
\end{figure}

We now compare the flow patterns of the DSO and DUE states, 
and discuss their relationships.
\cref{fig:DSO-DUE_FlowPtn} illustrates the arrival flow pattern at the destination in the DSO state (upper figure) and DUE state (lower figure).
Both arrival flow patterns are similar in terms of arrival time windows of commuters and the only difference is the destination arrival flow rates at every moment.
	Specifically, the arrival time windows of $(i,k)$-commuters in the DSO and DUE states are equal.
  Therefore, the sorting properties hold for both the DSO and DUE states.
However, the destination arrival flow rates of the DUE state differ from those of the DSO state because of the congestion effect of $\E{\dot{\sigma}}_{i}(t)$.
According to the definition of $\sigma_{i}(t)$, the congestion effect $\E{\dot{\sigma}}_{i}(t)$ is determined by the queueing delay pattern, and we find that the queueing delay pattern depends on the schedule delay cost function according to \cref{eq:DSO_sol_q,eq:DUE_final_sol_q}.
  
\par
We discuss the pricing patterns in the DSO state and queueing patterns in the DUE state.
We find that Pareto improvement can be achieved if the road manager imposes dynamic pricing that mimics the queueing delay pattern in the DUE state.
Two significant properties of \cref{pro:DUE_Solution} confirm this fact: 1) the DSO pricing patterns are the same as the DUE queuing delay patterns at each bottleneck; and 2) the DSO trip cost of each commuter (schedule delay costs + congestion prices) is equal to the DUE trip cost (schedule delay costs + queuing delay costs).

  Based on these properties, we find that the road manager's cost can be decreased by imposing the pricing equal to the queuing delay because the road manager gains revenue equal to the total queuing delay cost. 
  In other words, the road manager's benefit improves without increasing anyone's equilibrium commuting cost, thereby achieving the Pareto improvement.
  This is summarized in the following theorem.
  
\begin{thm}[Pareto improvement]
    Suppose that the QRP condition \eqref{eq:QRP_Corridor} holds.
		If the road manager imposes the dynamic pricing equal to the queuing delay $\{ \E{w}_{i}(t) \}_{t \in \ClT}$ at all bottlenecks $i \in \ClN$, the road manager's benefit improves without increasing anyone's equilibrium commuting cost
		\label{thm:ParetoImp_all_BN}
\end{thm}
\noindent
These results do not contradict those of previous studies that proposed that the optimal pricing does not achieve Pareto improvement when heterogeneous commuters exist (e.g.,~\citet{Arnott1988-pl,Arnott1994-dy,Van_den_Berg2010-ty,Van_den_Berg2011-qb,Hall2021-xo}).
	This is because these studies assume commuters' heterogeneity in terms of the value of time (i.e., parameter $\alpha$ in \cref{eq:TripCost} differs for each group, unlike the model of this study.
	\cref{thm:ParetoImp_all_BN} proposes that such inequity does not occur when the commuters' heterogeneity is considered in terms of the value of schedule delay.

\begin{figure}[tbp]
	\center
  \includegraphics[clip, width=0.65\columnwidth]{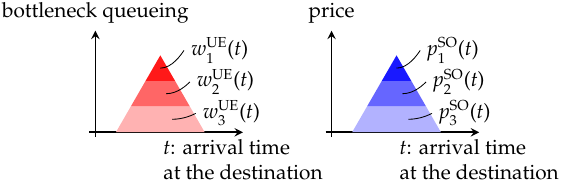}
	\caption{Simplified graphical images of bottleneck queueing and pricing patterns in a three-link network.}
	\label{fig:Queue_and_Price_Image}
\end{figure}

\begin{figure}[tbp]
	\center
  \includegraphics[clip, width=0.85\columnwidth]{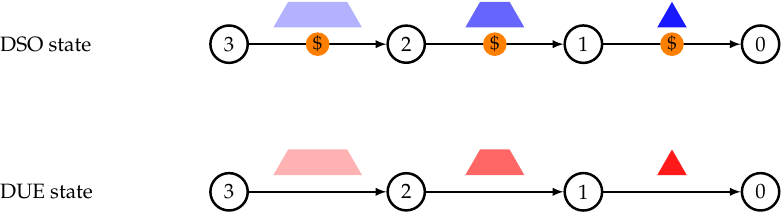} \caption{Example of the QRP in a three-link network.}
	\label{fig:QRP_3link_network}
\end{figure}

\subsection{Application of the QRP to the partial bottleneck pricing}
  Thus far, we have discussed that the QRP enables us to derive the DUE state by replacing the optimal pricing pattern with a bottleneck queueing delay pattern. 
	In this section, we focus on the case when pricing is introduced only for some bottlenecks and consider the equilibrium under the partial bottleneck pricing.
\par
	We consider the case in which the optimal pricing is 
	introduced to some bottlenecks $\ClN^{\mathrm{P}} \subseteq \ClN$ (hereafter referred to as partial bottleneck pricing: PBP).
	In this case, the trip cost of $(i,k)$-commuters arriving at the destination at time $t$ is calculated as follows:
\begin{align}
		&\PBP{C}_{i,k}(t) = 
		\sum_{j; j\leq i} \PBP{w}_{j}(t) 
		+ \left( \beta^{k} c(t) + \sum_{j; j \in \ClN^{\mathrm{P}}, j \leq i} 
		\S{p}_{j}(t) \right) + d_{i}
		&&\forall k \in \ClK, \quad 
			\forall i \in \ClN, \quad 
			\forall t \in \ClT,
\end{align}
where $\PBP{w}_{i}(t)$ represents the queueing delay at bottleneck $i$ for the commuters arriving at time $t$ at the destination.
The equilibrium state under the PBP can be obtained as the solution to the following LCP:
  \begin{align}
    \problemname{DUE-PBP}
  \\
	\mbox{Find} \quad
	& 
	\{\PBP{\Vtq}(t)\}_{t \in \ClT},
	\{\PBP{\Vtw}(t)\}_{t \in \ClT},
	\PBP{\Vtrho},
	\notag
	\\
	\mbox{such that} \quad
	& 
	0 \leq \PBP{q}_{i, k}(t) \perp
	\bracket{ 
		\PBP{C}_{i,k}(t) - \PBP{\rho}_{i, k}}
	  \geq 0
	 &&\forall k \in \ClK,
	 \quad \forall i \in \mathcal{N},
	 \quad \forall t \in \ClT,
	 \label{eq:DUE-PBP_DTC}
	 \\
	 &
	 0 \leq \PBP{w}_{i}(t) \perp 
	 \bracket{\mu_{i} \PBP{\dot{\sigma}}_{i}(t) - \sum_{j;j\geq i }  \PBP{q}_{j}(t) }
	 \geq  0
	 &&\forall i \in  \ClN \setminus \ClN^{\mathrm{P}},
	 \quad \forall t \in \ClT,
	 \label{eq:DUE-PBP_Queue}
	 \\
	 &0 \leq \PBP{\rho}_{i,k} \perp \bracket{\int_{\ClT}  \PBP{q}_{i, k}(t) \mathrm{d} t -  Q_{i, k}}
	 \geq 0
	 &&\forall k \in \ClK,
	 \quad \forall i \in \mathcal{N},
	 \label{eq:DUE-PBP_DemandCnsvCondition}
	 \\
	 &\PBP{\dot{\tau}}_{i}(t) > 0
	 &&\forall i \in \ClN, \quad \forall t \in \ClT.
	 \label{eq:DUE-PBP_ConsistencyCondition}
 \end{align}
	 The first condition \eqref{eq:DUE-PBP_DTC} represents the departure time choice condition.
	 The second condition \eqref{eq:DUE-PBP_Queue} represents the dynamics of the queueing at bottleneck $i \in \ClN \setminus \ClN^{\mathrm{P}}$.
	 The third condition \eqref{eq:DUE-PBP_DemandCnsvCondition} represents the demand conservation condition.
	 The final condition \eqref{eq:DUE-PBP_ConsistencyCondition} represents the consistency condition.
  \par
    Interestingly, if the pricing bottleneck set $\ClN^{\mathrm{P}}$ satisfies a certain condition, the queueing delay pattern at no-pricing bottleneck under the PBP equals that in the DUE state, i.e., $\PBP{w}_{i}(t)=\E{w}_{i}(t)$, $\forall i \in \ClN \setminus \ClN^{\mathrm{P}}$.
    Moreover, using this queueing delay pattern, we can construct the equilibrium flow pattern under the PBP, just like deriving the DUE solution. 
    To derive the solution in this way, we introduce the following assumption:
 \begin{asm}
      The pricing bottlenecks are located contiguously from the most upstream bottleneck $N$.
      That is, the pricing bottleneck set $\ClN^{\mathrm{P}} \subseteq \ClN$ can always be written as follows:
      \begin{align}
        \ClN^{\mathrm{P}} = \left\{ \PBP{i}, \PBP{i}+1, ..., N  \right\},
      \end{align}
      where $\PBP{i} \in \ClN$ is the most downstream pricing bottleneck.
   \label{asm:DownstreamSubset}
 \end{asm}
 \noindent
 In situations where this assumption holds, the queueing pattern under the PBP satisfies the consistency condition \eqref{eq:DUE-PBP_ConsistencyCondition}, which is the condition to guarantee a physically feasible flow in the Lagrangian coordinate approach as mentioned in \cref{sec:Model_settings}.
 Indeed, $\PBP{\dot{\tau}}_{i}(t)$ is calculated as follows:
 \begin{align}
  \PBP{\dot{\tau}}_{i}(t) &= 1  - \sum_{j;j\leq i } \PBP{\dot{w}}_{j}
  = 
  \begin{dcases}
    \E{\dot{\tau}}_{i}(t)
    \quad &\mathrm{if}\quad 
    i \in \ClN \setminus \ClN^{\mathrm{P}}
    \\
    \E{\dot{\tau}}_{(\PBP{i}-1)}(t)
    \quad &\mathrm{if}\quad
    i \in \ClN^{\mathrm{P}}
  \end{dcases}
  &&\forall i \in \ClN \setminus \{ 1 \}, 
  \quad \forall t \in \ClT,
 \end{align}
 because the PBP eliminates the queue at the pricing bottleneck $i \in \ClN^{\mathrm{P}}$ and maintains the queue at the unpriced bottleneck $i \in \ClN \setminus \ClN^{\mathrm{P}}$.
 Because $\E{\dot{\tau}}_{i}(t) > 0$, $\forall i \in \ClN$, $\forall t \in \ClT$, we can confirm that the consistency condition \eqref{eq:DUE-PBP_ConsistencyCondition} holds under \cref{asm:DownstreamSubset}. 
  Thus, this assumption is the sufficient condition for the queueing delay pattern $\E{w}_{i}(t)$, $\forall i \in \ClN \setminus \ClN^{\mathrm{P}}$ not to contradict the equilibrium condition under the PBP.\footnote{
    As a counterexample, we consider the case where \cref{asm:DownstreamSubset} does not hold.
      For example, consider the case where $\ClN=\{1,2,3\}$ and $\ClN^{\mathrm{P}} = \{ 2 \}$, it is not guaranteed that $\dot{\PBP{\tau}}_{3}(t)= 1- \E{\dot{w}}_{3}(t) - \E{\dot{w}_{1}}(t)$ satisfy consistency condition $\dot{\PBP{\tau}}_{3}(t)>0$ for all $t \in \ClT$.
}
 \par
   If we accept \cref{asm:DownstreamSubset}, we obtain the solutions to [DUE-PBP] using the QRP:
\begin{lem}
    Suppose that the QRP condition \eqref{eq:QRP_Corridor} and \cref{asm:DownstreamSubset} hold.
    Then, the following $\{\PBP{\Vtq}(t) \}_{t \in \ClT}$, $\{ \PBP{\Vtw}(t) \}_{t \in \ClT}$ and $\PBP{\Vtrho}$ are the solutions to [DUE-PBP]:
	\begin{align}
			&\PBP{q}_{i, k}(t) 
			= 
			\begin{dcases}
				\mu_{i} \PBP{\dot{\sigma}}_{i}(t) - \mu_{(i+1)} \PBP{\dot{\sigma}}_{(i+1)}(t)
				\quad &\mathrm{if} \quad t \in \E{\ClT}_{i, k} \setminus \E{\ClT}_{i, (k-1)}
        \\
        0
          &\mathrm{otherwise}
			\end{dcases}
			&&\forall i \in \ClN,
			\label{eq:PBP_q}
			\\
		&\PBP{\rho}_{i, k} = \E{\rho}_{i, k} 
		&&\forall i \in \ClN,
		\label{eq:PBP_rho}
		\\
		&\PBP{w}_{i}(t) 
		= 
		\begin{dcases}
			\E{w}_{i}(t)
			\quad &\mathrm{if} \quad i \in \ClN \setminus \ClN^{\mathrm{P}}
			\\
			0
			\quad &\mathrm{if} \quad i \in \ClN^{\mathrm{P}}
		\end{dcases}
		&&\forall t \in \ClT.
		\label{eq:PBP_w}
	\end{align}
	\label{lem:PBP_Solution}
\end{lem}
\noindent
\cref{lem:PBP_Solution} can be proven using a similar approach to the proof of \cref{thm:QRP_Corridor}.
Specifically, when queues at the bottleneck with pricing are replaced by zero, and the queues at the bottleneck without pricing are replaced by $p$, we show that these queueing delay patterns satisfy the DUE condition under the PBP.
\par
\begin{figure}[tbp]
	\center
  \includegraphics[clip, width=0.65\columnwidth]{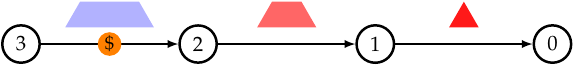}
	\caption{Equilibrium under the partial bottleneck pricing (PBP).}
	\label{fig:partial_bottleneck_pricing_1}
\end{figure}

Let us consider the queueing delay patterns that occur under the PBP.
	For convenience, we introduce \cref{fig:Queue_and_Price_Image,fig:QRP_3link_network}, which illustrate a simple example of the QRP in a three-link network.
	\cref{fig:Queue_and_Price_Image} depicts the pricing and queueing patterns in a three-link network.
	The approximate shapes of both patterns are represented by triangles and trapezoids.
	Using these graphic representations, 
	\cref{fig:QRP_3link_network} shows the QRP in which the DSO pricing patterns are the same as the DUE queuing delay patterns at each bottleneck.
    \cref{fig:partial_bottleneck_pricing_1} depicts an equilibrium state under the pricing at Bottleneck $3$, i.e., $\ClN^{\mathrm{P}} = \{ 3 \}$.
    	Because the road manager imposes tolling $\S{p}_{3}(t)$ at Bottleneck $3$, the queueing at Bottleneck $3$ is completely eliminated.
	Conversely, the bottleneck queueing is formed at Bottlenecks $1$ and $2$, thus, maintaining the same equilibrium commuting cost.
	Consequently, the bottleneck queueing at Bottlenecks $1$ and $2$ is shown in \cref{fig:partial_bottleneck_pricing_1}.
	As shown in \cref{lem:PBP_Solution}, this equilibrium state has two essential features: (i) the queueing delay pattern at bottleneck $i \in \ClN \setminus \ClN^{\mathrm{P}}$ is equal to that in the DUE state, and (ii) the equilibrium commuting cost of each commuter is equal to that in the DUE (DSO) state.

By comparing the solutions to the [DSO], [DUE], and [DUE-PBP], we obtain the following theorem:
\begin{thm}
		Suppose that the QRP condition \eqref{eq:QRP_Corridor} and \cref{asm:DownstreamSubset} hold.
		\par \noindent
	(a) The equilibrium commuting cost of $(i,k)$-commuters under the 
	PBP is equal to that in the DSO and DUE states:
	\begin{align}
    &\PBP{\rho}_{i,k} = 
		\S{\rho}_{i,k} = \E{\rho}_{i,k}
		&&\forall k \in \ClK, \quad \forall i \in \ClN.
		\label{eq:PBP_EquilibriumCost_are_equal}
	\end{align}
	\par \noindent
	(b)	The total system costs satisfy the following inequality:
	\begin{align}
		&\S{Z} \leq \PBP{Z}(\ClN^{\mathrm{P}}) \leq \E{Z},
		\label{eq:PBP_TotalCost_DSO<PBP}
	\end{align}
    where $\S{Z}$ and $\E{Z}$ are the total system costs of the DSO and DUE states, respectively.
      The equality in the LHS of \cref{eq:PBP_TotalCost_DSO<PBP} holds if and only if $\ClN^{\mathrm{P}}=\ClN$ and the equality in the RHS of \cref{eq:PBP_TotalCost_DSO<PBP} holds if and only if $\ClN^{\mathrm{P}}=\emptyset$.
	\label{thm:PBP_Pareto}
\end{thm}
\noindent
\cref{thm:PBP_Pareto} implies that we can improve the system cost by partially implementing dynamic pricing on bottlenecks without considering the influence of such an implementation on the other bottlenecks.
In addition, the road manager gains revenue equal to the sum of queuing delay costs at the bottlenecks with pricing.
In other words, the partial bottleneck pricing can achieve the Pareto improvement under \cref{asm:DownstreamSubset} and the QRP condition. 
Therefore, because of \cref{thm:ParetoImp_all_BN,thm:PBP_Pareto}, we can interpret the QRP condition as a sufficient condition for a Pareto improvement by full/partial optimal pricing.
A similar result was shown by~\citet{Fu2022-nl}, assuming homogeneous commuters; however, the remarkable difference is that the pricing bottleneck set condition (\cref{asm:DownstreamSubset}) is required.

\begin{figure}[tbp]
	\center
  \includegraphics[clip, width=0.8\columnwidth]{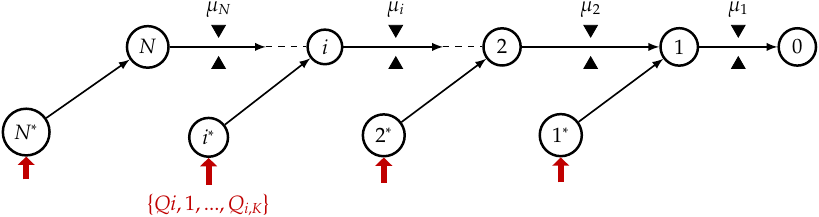}
	  \caption{Corridor network with on-ramp links.}
    \label{fig:RampMetering_Network}
\end{figure}

\section{Applications of the QRP to on-ramp-based policies}\label{sec:Application}
The QRP and PBP analyses results represent a replaceability of commuters' costs in the departure time choice problem under several conditions.
This replaceability means that as long as each commuter's commuting cost is equal to those in the DUE and DSO states, a feasible flow pattern exists regardless of whether the costs experienced at the links are queuing or price costs.
Such link cost replaceability leads us to expect the same kind of replaceability for the node-based (on-ramp-based) cost. 
For example, the queueing at the on-ramp under the on-ramp metering may coincide with the pricing pattern under the on-ramp pricing. 
If such replaceability (equivalence) also holds, we can construct a systematic approach for analyzing the on-ramp-based and bottleneck-based policies based on the QRP.
\par
	This section considers the two on-ramp-based policies: on-ramp metering and on-ramp pricing to gain an insight into the theoretical relationship between the on-ramp-based and bottleneck-based policies.
    We demonstrate that the QRP enables us to obtain equilibrium solutions under the on-ramp-based policies. 
	\cref{subsec:Full_implementation} analyzes the case in which metering/pricing is implemented at all on-ramps (hereafter referred to as full implementation) and derives the equilibrium solution using the QRP.
    \cref{subsec:Partial_implementation} analyzes the case when metering/pricing is only implemented at partial on-ramps (hereafter referred to as partial implementation) and derives the equilibrium solution to this case.
	\cref{subsec:Comparing} discusses the theoretical properties of these policies by comparing these equilibrium solutions in the DSO and DUE solutions.
		In the remainder of this section, we assume that the QRP schedule delay cost function satisfies \cref{eq:QRP_Corridor}.
\par
We consider the extended network with on-ramp links and nodes.
As shown in \cref{fig:RampMetering_Network}, 
we add on-ramp link $i^{\ast}$ representing the access link that connects nodes $i$ and $i^{\ast}$.
Free-flow travel times of on-ramp links are assumed to be zero, and the set of on-ramp links is defined by $\ClN^{\ast}$.
The arrival flow variables are defined by $q_{i^{\ast}}(t)$ using $i^{\ast}$ 
instead of $q_{i}(t)$.
We define $w_{i^{\ast}}(t)$ as the queueing delay at on-ramp $i^{\ast}$ created by the on-ramp metering for the commuters whose destination arrival time is $t$.

\subsection{Full implementation policies}\label{subsec:Full_implementation}

\subsubsection{Full on-ramp metering}
On-ramp metering (RM) is one of the most common methods used for managing freeway traffic.
The primary strategy of ramp metering 
is to maintain the freeway traffic by adjusting the inflow rate from the on-ramp.
This study considers ramp metering, which aims to eliminate queues on the freeway.
  For all on-ramps $i^{\ast} \in \ClN^{\ast}$, we set the metering rate to $\mu_{i^{\ast}} = \overline{\mu}_{i}$ because $\overline{\mu}_{i}$ is the maximum value without creating a queue at the downstream bottleneck $i$ on the freeway.
\par
We consider the equilibrium under the on-ramp metering for all ramps (full on-ramp metering).
The associated equilibrium problem is the following LCP:
\begin{align}
  \problemname{DUE-RM}
  \\
\mbox{Find} \quad
 &
 \{\Vtq(t)\}_{t \in \ClT},
 \{\Vtw(t)\}_{t \in \ClT},
    \Vtrho
\notag
 \\
 \mbox{such that} \quad
 &
 0 \leq q_{i^{\ast}, k}(t) \perp
 \bracket{ w_{i^{\ast}}(t) +  \beta^{k} c(t) + d_{i}
 - \rho_{i^{\ast}, k} }
 \geq  0
 &&\forall k \in \ClK,
 \quad \forall i \in \mathcal{N},
 \quad \forall t \in \ClT,
 \label{eq:RM-DTC}
 \\
 &
 0 \leq w_{i^{\ast}}(t) \perp 
 \bracket{ \overline{\mu}_{i} - q_{i^{\ast}}(t) }
\geq 0
 &&\forall i^{\ast} \in \ClN^{\ast},
	\quad \forall t \in \ClT,
\label{eq:RM-Queue}
	\\
& 0 \leq \rho_{i^{\ast},k} \perp
\bracket{ \int_{\ClT}  q_{i^{\ast}, k}(t) \mathrm{d} t
- Q_{i, k}}
\geq 0
  &&\forall k \in \ClK,
	\quad \forall i^{\ast} \in \mathcal{N}^{\ast},
\label{eq:RM-DemandCnsv}
	\\
	&\dot{\tau}_{i^{\ast}}(t)  > 1
  &&\forall i^{\ast} \in \ClN^{\ast}.
\label{eq:RM-Consistency}
\end{align}
	Condition \eqref{eq:RM-Queue} represents the queueing dynamics at each on-ramp $i^{\ast} \in \ClN^{\ast}$. 
\par
\begin{figure}[tbp]
	\center
  \includegraphics[clip, width=0.70\columnwidth]{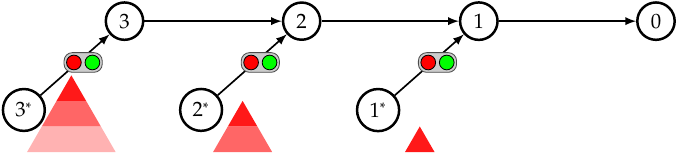}
	\caption{Equilibrium under the on-ramp metering (RM).}
	\label{fig:ramp_metering}
\end{figure}

  The problem [DUE-RM] has a similar mathematical structure to the problem [DUE], particularly, for the departure time choice condition of [DUE-RM] and [DUE] (\cref{eq:RM-DTC,eq:DUE-DTC}), if $w_{i^{\ast}}(t)$ is regarded as $\sum_{j;j\leq i} \E{w}_{j}(t)$, then the form is exactly the same.
    This fact leads us to expect each commuter's commuting costs of [DUE-RM] to be equal to those of [DUE].
    This expectation is actually true.
    We confirm this by constructing the solution to [DUE-RM] using the equilibrium condition of [DUE-RM] and the solution to [DUE]. 
  The solution to [DUE-RM] is shown as follows:
\begin{lem}
 The following $\{\RM{\Vtq}(t)\}_{t \in \ClT}$, $\{\RM{\Vtw}(t)\}_{t \in \ClT}$, and $\RM{\Vtrho}$ are equilibrium solutions to [DUE-RM]:
	\begin{align}
		&\RM{q}_{i^{\ast}, k}(t)
		=
		\begin{dcases}
      \overline{\mu}_{i} 
        &\mathrm{if}\quad
       t \in \E{\ClT}_{i, k} \setminus \E{\ClT}_{i, (k-1)}
      \\
      0
        &\mathrm{otherwise}
    \end{dcases}
		&&\forall k \in \ClK, \quad \forall i \in \ClN, \quad \forall t \in \ClT.
		\\
		&\RM{w}_{i^{\ast}}(t) = \sum_{j;j\leq i} \E{w}_{j}(t)
		&&\forall i^{\ast} \in \ClN^{\ast}, \quad \forall t \in \ClT,
		\\
		&\RM{\rho}_{i, k} = \E{\rho}_{i, k}
		&&\forall k \in \ClK, \quad \forall i \in \ClN.
	\end{align}
	\label{lem:RM_solution}
\end{lem}
  \cref{fig:ramp_metering} depicts the solution to [DUE-RM] in the case of a three-link network.
  We find that the queueing delay at each on-ramp $\RM{w}_{i^{\ast}}(t) $ is equal to the sum of queueing delays at its downstream bottlenecks in the DUE state $\sum_{j;j\leq i} \E{w}_{j}(t)$.

\begin{figure}[tbp]
	\center
  \includegraphics[clip, width=0.70\columnwidth]{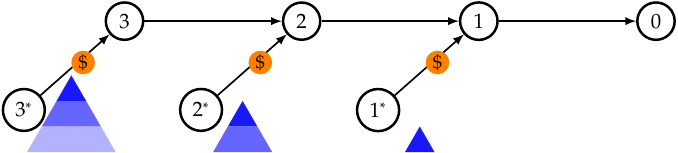}
	\caption{Equilibrium under the on-ramp pricing (RP).}
	\label{fig:ramp_pricing}
\end{figure}

\subsubsection{Full on-ramp pricing}
We consider the equilibrium under the on-ramp pricing.
In this case, the road manager imposes the toll at 
all on-ramps. 
We assume that the road manager tolled $\RP{p}_{i^{\ast}}(t) = \sum_{j;j\leq i} \S{p}_{j}(t)$ at on-ramp  $i^{\ast}$.
The associated equilibrium problem is the following LCP:
\begin{align}
  \problemname{DUE-RP}
  \\
 \mbox{Find} \quad
 &
 \{ \Vtq(t)\}_{t \in \ClT},
	\{ \Vtw(t)\}_{t \in \ClT},
	\Vtrho
\notag
 \\
 \mbox{such that} \quad
 &
 0 \leq q_{i^{\ast}, k}(t) \perp 
 \bracket{
	\RP{p}_{i^{\ast}}(t) + w_{i^{\ast}}(t) +  \beta^{k} c(t)
	+ d_{i} - \rho_{i^{\ast}, k} }
	\geq  0
	&&\forall k \in \ClK,
	\quad \forall i \in \mathcal{N},
	\quad \forall t \in \ClT,
	\label{eq:RP-DTC}
	\\
	&
	0 \leq w_{i^{\ast}}(t) \perp
	\bracket{ \overline{\mu}_{i} - \RP{q}_{i^{\ast}}(t)}
	\geq 0
	&&\forall i^{\ast} \in \ClN^{\ast},
	\quad \forall t \in \ClT,
	\label{eq:RP-Queue}
	\\
	& 0 \leq \rho_{i^{\ast},k} \perp
	\bracket{ \int_{\ClT}  q_{i^{\ast}, k}(t) \mathrm{d} t
	- Q_{i, k}}
	\geq 0
  &&\forall k \in \ClK, 
  \quad \forall i^{\ast} \in \ClN^{\ast},
	\label{eq:RP-DemandCnsv}
	\\
	&\dot{\tau}_{i^{\ast}}(t)  > 1
  &&\forall i^{\ast} \in \ClN^{\ast}.
	\label{eq:RP-Consistency}
\end{align}
	Condition \eqref{eq:RP-Queue} represents the queueing dynamics at each bottleneck $i \in \ClN$. 
	The queues at all ramps are completely eliminated by the on-ramp pricing ($\RP{w}_{i^{\ast}}(t) = 0$, $\forall i^{\ast} \in \ClN^{\ast}$, $\forall t \in \ClT$).
\par
  Using an approach similar to that of [DUE-RM], we obtain the solution to [DUE-RP] as follows:
\begin{lem}
    The following $\{\RP{\Vtq}(t)\}_{t \in \ClT}$, $\{\RP{\Vtw}(t)\}_{t \in \ClT}$, and $\RP{\Vtrho}$ are equilibrium solutions to [DUE-RP]:
	\begin{align}
		&\RP{q}_{i^{\ast}, k}(t)
		=
		\begin{dcases}
      \overline{\mu}_{i} 
        &\mathrm{if}\quad
      t \in \E{\ClT}_{i, k} \setminus \E{\ClT}_{i, (k-1)}
      \\
      0
        &\mathrm{otherwise}
    \end{dcases}
		&&\forall k \in \ClK, \quad \forall i \in \ClN, \quad \forall t \in \ClT.
		\\
		&\RP{w}_{i^{\ast}}(t) = 0
		&&\forall i^{\ast} \in \ClN^{\ast}, \quad \forall t \in \ClT,
		\\
		&\RP{\rho}_{i, k} = \E{\rho}_{i, k} 
		&&\forall k \in \ClK, \quad \forall i \in \ClN,
	\end{align}
	\label{lem:RP_solution}
\end{lem}
\noindent
\cref{fig:ramp_pricing} depicts the equilibrium state under the on-ramp pricing.
We see that the queueing is equal to the on-ramp pricing at each on-ramp.

\subsection{Partial implementation policies}\label{subsec:Partial_implementation}

\subsubsection{Partial on-ramp metering}
We consider the partial on-ramp metering and define the set of links with the on-ramp meter as $\ClN^{\ast\mathrm{P}} \subseteq \ClN^{\ast}$.
The associated equilibrium problem is the following LCP:
\begin{align}
  \problemname{DUE-PRM}
  \\
	\mbox{Find} \quad
 &
   \{ \Vtq(t) \}_{t \in \ClT}
	 \{ w_{i}(t) \}_{i \in \ClN, t \in \ClT},
	 \{ w_{i^{\ast}}(t) \}_{i^{\ast} \in \ClN^{\ast}, t \in \ClT},
	 \Vtrho,
\notag
 \\
 \mbox{such that} \quad
	&
	0 \leq q_{i^{\ast}, k}(t) \perp
	\bracket{ w_{i^{\ast}}(t) + \sum_{j;j<i} w_{j}(t) +  \beta^{k} c(t)
	+ d_{i}	- \rho_{i^{\ast}, k} }
	\geq  0
&&
  \forall i^{\ast} \in \ClN^{\ast},
  \forall k \in \ClK,
		\quad \forall t \in \ClT,
	\label{eq:PRM-DTC}
	\\
	&
	0 \leq w_{i}(t) \perp
	\bracket{\mu_{i} \dot{\sigma}_{i}(t) - \sum_{j^{\ast};j^{\ast} \geq i^{\ast} } q_{j^{\ast}}(t)}
 \geq 0
	&&\forall i \in \ClN \setminus \ClN^{\ast\mathrm{P}},
	 \quad \forall t \in \ClT,
	\label{eq:PRM-Queue-1}
	\\
	&
	0 \leq w_{i^{\ast}}(t) \perp
	\bracket{\overline{\mu}_{i} \dot{\sigma}_{i^{\ast}}(t) - q_{i^{\ast}}(t)}
 \geq 0
 &&\forall i^{\ast} \in \ClN^{\ast\mathrm{P}},
	 \quad \forall t \in \ClT,
	\label{eq:PRM-Queue-2}
	\\
	&0 \leq \rho_{i,k} \perp 
	\bracket{ \int_{\ClT} q_{i^{\ast}, k}(t) \mathrm{d} t - Q_{i, k} }
	\geq 0
	&&\forall k \in \ClK,
	\quad \forall i \in \mathcal{N},
	\label{eq:PRM-DemandCnsv}
	\\
	&\dot{\tau}_{j}(t) > 0
	&&\forall j \in \ClN \cup \ClN^{\ast},
	\quad \forall t \in \ClT.
	\label{eq:PRM-Consistency}
\end{align}
	The first condition \eqref{eq:PRM-DTC} represents the departure time choice condition of each $(i^{\ast}, k)$-commuters.
	Commuters who depart from origin $i^{\ast} \in \ClN^{\ast\mathrm{P}}$ are constrained by the on-ramp metering.
	Conditions \eqref{eq:PRM-Queue-1} and \eqref{eq:PRM-Queue-2} represent the queueing dynamics at each bottleneck $i \in \ClN \setminus \ClN^{\ast\mathrm{P}}$ and each on-ramp $i \in \ClN^{\ast\mathrm{P}}$, respectively.
\par
\begin{figure}[tbp]
  \begin{minipage}[b]{0.99\columnwidth}
	\center
  \includegraphics[clip, width=0.70\columnwidth]{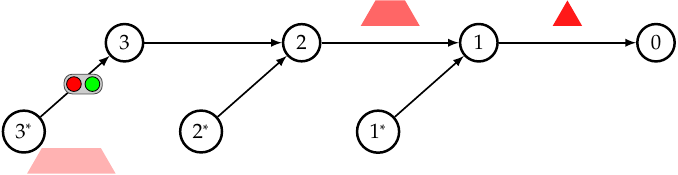}
  \subcaption{Case of $\ClN^{\ast\mathrm{P}} = \{ 3 \}$}
  \label{fig:partial_metering_1}
\end{minipage}
\begin{minipage}[b]{0.99\columnwidth}
	\center
  \includegraphics[clip, width=0.70\columnwidth]{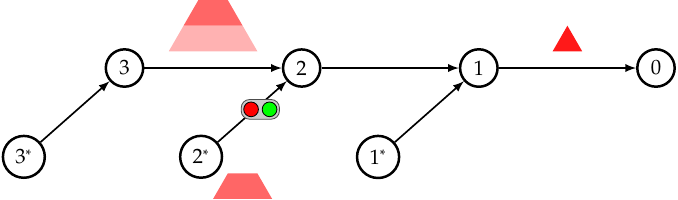}
	\subcaption{Case of $\ClN^{\ast\mathrm{P}} = \{ 2 \}$}
	\label{fig:partial_metering_2}
\end{minipage}
\caption{Equilibrium under the partial on-ramp metering (PRM).}
\label{fig:PRM}
\end{figure}

  Although the [DUE-PRM] appears to be a complex equilibrium problem at first glance, we obtain the solution by applying the methods discussed previously.
The basic approach for obtaining the solution is the same as in the case of partial bottleneck pricing. 
Specifically, the approach is based on the conjecture that the equilibrium commuting cost of each commuter is equal to that in the DUE (DSO) state. 
Based on this conjecture, we derive the solution to [DUE-PRM] using the equilibrium condition of [DUE-PRM] and the DUE solution. 
We subsequently confirm that the conjecture is true.
Consequently, we analytically obtain the solutions to [DUE-PRM]:
\begin{lem}
   The following 
   $\{ \PRM{\Vtq}(t) \}_{t \in \ClT}$
	 $\{ \PRM{w}_{i}(t) \}_{i \in \ClN, t \in \ClT}$,
	 $\{ \PRM{w}_{i^{\ast}}(t) \}_{i^{\ast} \in \ClN^{\ast}, t \in \ClT}$, and $\PRM{\Vtrho}$
	 are the equilibrium solution to [DUE-PRM]:
	\begin{align}
		&\PRM{\rho}_{i, k} = \E{\rho}_{i, k} 
		&&\forall i \in \ClN,
		\\
		&\PRM{w}_{i}(t) 
		= 
		\begin{dcases}
			0
			\quad &\mathrm{if} \quad i \in \ClN^{\ast\mathrm{P}}
			\\
			\sum_{j \in  \ClN^{\ast\mathrm{P}}; j \leq i } \E{w}_{j}(t)
			\quad &\mathrm{if} \quad i \in\ClN \setminus \ClN^{\ast\mathrm{P}}
		\end{dcases}
		&&\forall t \in \ClT,
		\\
		&\PRM{w}_{i\ast}(t) 
		= 
		\begin{dcases}
			\sum_{j \in \ClN^{\ast} \setminus \ClN^{\ast\mathrm{P}} ; j \leq i} \E{w}_{j}(t)
			\quad &\mathrm{if} \quad i \in \ClN^{\ast\mathrm{P}}
			\\
			0
			\quad &\mathrm{if} \quad i \in \ClN^{\ast} \setminus \ClN^{\ast\mathrm{P}}
		\end{dcases}
		&&\forall t \in \ClT,
		\\
		&\PRM{q}_{i^{\ast}, k}(t) 
			= 
			\begin{dcases}
				\mu_{i} \PRM{\dot{\sigma}}_{i}(t) - \mu_{(i+1)} \PRM{\dot{\sigma}}_{(i+1)}(t)
				\quad &\mathrm{if}\quad t \in \E{\ClT}_{i, k}
				\\
				0
				\quad &\mathrm{if}\quad t \notin \E{\ClT}_{i, k}
			\end{dcases}
			&&\forall i \in \ClN.
	\end{align}
	\label{lem:PRM_Solution}
\end{lem}
\cref{fig:partial_metering_1} depicts the equilibrium state under the on-ramp metering at Link $3$, i.e., $\ClN^{\ast\mathrm{P}} = \{ 3 \}$.
The on-ramp queueing delay pattern at Link $3$ is equal to 
the pattern under all on-ramp metering cases, and there is no queueing at Bottleneck $3$.
Moreover, the bottleneck queues are formed at Bottlenecks $1$ and $2$, maintaining the same equilibrium commuting cost.
The bottleneck queues at Bottlenecks $1$ and $2$ are shown in \cref{fig:partial_metering_1}.
Any combination of ramps with metering is acceptable.
For example, \cref{fig:partial_metering_2} depicts the equilibrium state under the on-ramp metering at Link $2$, i.e., $\ClN^{\ast\mathrm{P}} = \{ 2 \}$.

\subsubsection{Partial on-ramp pricing}

\begin{figure}[tbp]
  \begin{minipage}[b]{0.99\columnwidth}
	\center
  \includegraphics[clip, width=0.70\columnwidth]{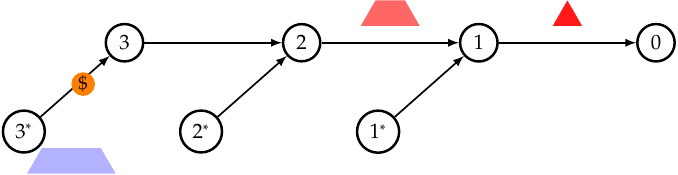}
	\subcaption{Case of $\ClN^{\ast\mathrm{P}} = \{ 3 \}$}
	\label{fig:partial_on-ramp_pricing_1}
\end{minipage}
\begin{minipage}[b]{0.99\columnwidth}
	\center
  \includegraphics[clip, width=0.70\columnwidth]{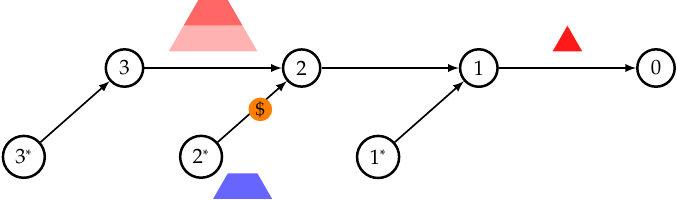}
  \subcaption{Case of $\ClN^{\ast\mathrm{P}} = \{ 2 \}$}
  \label{fig:partial_on-ramp_pricing_2}
\end{minipage}
\caption{Equilibrium under the partial on-ramp pricing (PRP).}
\label{fig:PRP}
\end{figure}
We now consider the partial on-ramp pricing.
Let $\ClN^{\ast\mathrm{P}}$ be a set of links with pricing.
Assuming that the road manager tolled $\RP{p}_{i^{\ast}}(t)$ at on-ramp $ i^{\ast} \in \ClN^{\ast\mathrm{P}}$, 
the associated equilibrium problem is the following LCP:
	\begin{align}
    \problemname{DUE-PRP}
    \\
		\mbox{Find} \quad
    &
    \{ \Vtq(t) \}_{t \in \ClT}
    \{  \Vtw(t) \}_{t \in \ClT},
    \Vtrho,
		\notag
		\\
		\mbox{such that} \quad
		&
      0 \leq q_{i^{\ast}, k}(t) \perp
		\bracket{\sum_{j; j\leq i} w_{j}(t) 
		 + \beta^{k} c(t) + \RP{p}_{i^{\ast}}(t)
		 + d_{i} - \rho_{i, k} }
		 \geq 0
		 &&\forall k \in \ClK,
		 \quad \forall i^{\ast} \in \ClN^{\ast\mathrm{P}},
		 \quad \forall t \in \ClT,
		 \label{eq:PRP-DTC}
		\\
		 &
		 0 \leq w_{i}(t) \perp 
		 \bracket{ \mu_{i} \dot{\sigma}_{i}(t) - \sum_{j;j\geq i } q_{j}(t) }
		 \geq  0
		 &&\forall i \in \mathcal{N},
		 \quad \forall t \in \ClT,
		 \label{eq:PRP-Queue}
		 \\
		 &0 \leq \rho_{i,k} \perp 
		 \bracket{ \int_{\ClT}  q_{i, k}(t) \mathrm{d} t - Q_{i, k}}
		 \geq 0
		 &&\forall k \in \ClK,
		 \quad \forall i \in \mathcal{N},
		 \label{eq:PRP-DemandCnsv}
		 \\
     &\dot{\tau}_{j}(t) > 0
     &&\forall j \in \ClN \cup \ClN^{\ast},
     \quad \forall t \in \ClT.
		 \label{eq:PRP-Consistency}
\end{align}
	Condition \eqref{eq:PRP-DTC} represents the departure time choice condition of each commuter. 
\par
    We also obtain the equilibrium solution to [DUE-PRP] using an approach similar to that used in the case of [DUE-PRM].
  Consequently, we obtain the following equilibrium solution:
\begin{lem}
    The following  
    $\{ \PRP{\Vtq}(t) \}_{t \in \ClT}$
    $\{  \PRP{\Vtw}(t) \}_{t \in \ClT}$,
    and $\PRP{\Vtrho}$ are the equilibrium solutions to [DUE-PRP]:
	\begin{align}
    &\PRP{q}_{i^{\ast}, k}(t) 
    = 
    \begin{dcases}
      \mu_{i} \PBP{\dot{\sigma}}_{i}(t) - \mu_{(i+1)} \PBP{\dot{\sigma}}_{(i+1)}(t)
      \quad &\mathrm{if}\quad t \in \E{\ClT}_{i, k} \setminus \E{\ClT}_{i, (k-1)}      
      \\
      0
      \quad &\mathrm{otherwise}
    \end{dcases}
    &&\forall i \in \ClN,
		\\
		&\PRP{w}_{i}(t) 
		= 
		\begin{dcases}
			0
			\quad &\mathrm{if} \quad i \in \ClN^{\ast\mathrm{P}}
			\\
			\sum_{j \in  \ClN^{\ast\mathrm{P}}; j \leq i } \E{w}_{j}(t)
			\quad &\mathrm{if} \quad i \in\ClN \setminus \ClN^{\ast\mathrm{P}}
		\end{dcases}
		&&\forall t \in \ClT,
    \\
    &\PRP{\rho}_{i, k} = \E{\rho}_{i, k} 
    &&\forall i \in \ClN.
  \end{align}
	\label{lem:PRP_Solution}
\end{lem}
\cref{fig:partial_on-ramp_pricing_1} and \cref{fig:partial_on-ramp_pricing_2} depict the solution to [DUE-PRP] in the case of $\ClN^{\ast\mathrm{P}} = \{ 3 \}$ and $\ClN^{\ast\mathrm{P}} = \{ 2 \}$, respectively.

\subsection{Comparison of policies}\label{subsec:Comparing}
  We find several insights into welfare analysis by comparing the equilibrium solutions under various policies that have been analytically derived thus far.
  First, the equilibrium commuting cost of $(i,k)$-commuters under the full implementation policies is equal to that of the DUE and DSO states.
		Based on this fact, the total cost under full implementation policies can be calculated.
		Therefore, we obtain the following theorem:
\begin{thm}
		Suppose that the QRP condition \eqref{eq:QRP_Corridor} holds.
	\par \noindent
	(a) The equilibrium commuting cost of $(i,k)$-commuters under the RP and RM are equal to that of the DUE and DSO states:
	\begin{align}
		\S{\rho}_{i,k} = \E{\rho}_{i,k} 
		= \RP{\rho}_{i,k}
		= \RM{\rho}_{i,k}.
		\label{eq:EquilibriumCost_are_equal_Full}
	\end{align}
	(b) The total system cost under the full implementation policies satisfies the following inequalities:
	\begin{align}
		\S{Z} = \RP{Z} < \RM{Z} = \E{Z},
		\label{eq:S-BP-RM-E_CostIneq}
	\end{align}
  where $\RM{Z}$ and $\RP{Z}$ are total system costs in equilibrium under the on-ramp metering and on-ramp pricing, respectively.
	\label{thm:TotalCost_DSO=RP<RM=DUE}
\end{thm}
\noindent
\cref{eq:EquilibriumCost_are_equal_Full} shows that the full implementation policies do not change each commuter's equilibrium commuting cost, like the bottleneck-based pricing.
Because this is true for all commuters, regardless of the origin or group, we see that introducing the full implementation policies neither creates those who lose nor gain.
This may improve the social acceptability of introducing these policies.
\par
\cref{thm:TotalCost_DSO=RP<RM=DUE} shows that if the QRP condition \eqref{eq:QRP_Corridor} holds, the bottleneck pricing at all bottlenecks simultaneously achieves the first best and Pareto improvement.
Conversely, on-ramp metering eliminates congestion 
on the freeway but creates queues on the on-ramp, 
resulting in total costs equal to those in the pure DUE state.
\par
We then focus on the partial implementation policies.
By comparing the equilibrium solutions under the partial implementation policies with the DUE and DSO states, we obtain the following theorem:
\begin{thm}
		Suppose that the QRP condition \eqref{eq:QRP_Corridor} holds.
	\par \noindent
	(a) The equilibrium commuting cost of $(i,k)$-commuters under the PRM and PRP is equal to that of DUE and DSO states:
	\begin{align}
		\S{\rho}_{i,k} = \E{\rho}_{i,k} 
		= \PRM{\rho}_{i,k} 
		= \PRP{\rho}_{i,k}.
		\label{eq:EquilibriumCost_are_equal_Partial_1}
	\end{align}
		Regardless of where the policy is implemented ($\ClN^{\mathrm{P}}$, $\ClN^{\ast\mathrm{P}}$), \cref{eq:EquilibriumCost_are_equal_Partial_1} holds.
	\par \noindent
	(b) The total system costs satisfy the following inequality:
	\begin{align}
		&\S{Z} \leq  \PRP{Z}(\ClN^{\ast\mathrm{P}}) \leq \E{Z},
		\label{eq:TotalCost_DSO<PRP<DUE}
		\\
		&\S{Z} < \PRM{Z}(\ClN^{\ast\mathrm{P}}) = \E{Z},
		\label{eq:TotalCost_DSO<PRM=DUE}
	\end{align}
  where $\PRM{Z}(\ClN^{\ast\mathrm{P}})$ and $\PRP{Z}(\ClN^{\ast\mathrm{P}})$ are the total system costs at equilibrium when on-ramp metering and pricing are implemented on parts of on-ramps $\ClN^{\ast\mathrm{P}}$, respectively.
	\label{thm:TotalCost_DSO<PRP<PRM=DUE}
\end{thm}

\begin{table}[tbp]
  \centering
  \caption{Comparison of the policies.}
  \label{tab:Policies_implications}
    \begin{tabular}{ccl||cc}
    \hline
		& &       & \multicolumn{1}{c}{Total System Cost} 
					& \multicolumn{1}{c}{Pareto improvement} \\
    \hline
    \hline
    \multirow{1}[4]{*}{Bottleneck-based} & \multirow{1}[4]{*}{Pricing} & Full  & $=$  DSO ($<$DUE) 
		& $\checkmark$ \\
    \cline{3-5} {}&         & Partial & $\geq$ DSO ($<$DUE) & $\checkmark$ \\
    \hline
    \multirow{3}[4]{*}{On-ramp-based} &\multirow{1}[4]{*}{Metering} & Full  & $=$ DUE  &  \\
    \cline{3-5} {}&         & Partial & $=$ DUE  &  \\
    \cline{2-5} 
    {}&\multirow{1}[4]{*}{Pricing} & Full  & $=$ DSO ($<$DUE) & $\checkmark$ \\
    \cline{3-5} {}&        & Partial (*) & $\geq$ DSO ($<$DUE) & $\checkmark$ \\
    \hline
	\end{tabular}%
	\\
	\rightline{(*) under \cref{asm:DownstreamSubset}}
\end{table}%

\cref{thm:TotalCost_DSO<PRP<PRM=DUE} provides significant insights.
First, the partial on-ramp pricing achieves the Pareto improvement because the total system cost decreases (i.e., the road manager's  benefits increase) without increasing anyone's equilibrium commuting costs (\cref{thm:TotalCost_DSO<PRP<PRM=DUE}).
This implies that partial pricing of on-ramps should be implemented whenever possible, even when only parts of a network can be priced, or pricing is constrained in other ways.
Second, we confirm that the partial on-ramp metering, similar to the full implementation case, causes a queue on the ramp while eliminating the queue at the direct downstream link on the freeway.
Consequently, the total system cost is equal to that of the DUE states.
\par
Results from this section are summarized in \cref{tab:Policies_implications}.
Comparing these results, we conclude that the on-ramp metering has no significant benefit.\footnote{
    In general, the on-ramp metering prevents the capacity drop on the freeway~\citep{Papageorgiou2002-jj}.
    We realize that if the effect of the capacity drop is considered in our model, the results may be updated.}
  The full pricing policies, whether bottleneck-based or on-ramp-based, can simultaneously achieve the system optimum and Pareto improvement.
  The on-ramp-based pricing may be more practical than the bottleneck-based pricing because it does not require any conditions for partial implementation.

\section{Conclusion}\label{sec:ConcludingRemarks}
This study investigated the DSO and DUE problems in a corridor network with multiple bottlenecks, considering the commuter heterogeneity with respect to schedule delay. 
Consequently, an analytical solution to the DSO and DUE problems was derived successfully using a systematic approach. 
We first formulated the DSO problem as an LP.
Subsequently, we derived the analytical solution to the DSO problem by applying the theory of optimal transport. 
Then, we demonstrated the QRP, that is, the DSO state without queues can be achieved by imposing a congestion price equal to the queuing delay in the DUE state at each bottleneck and 
at every moment under certain conditions of a schedule delay cost function. 
Finally, based on the DSO solution and QRP, 
we derived an analytical solution to the DUE problem, which was formulated as an LCP. 
Moreover, as an application of the QRP and analytical solutions, we demonstrated that the equilibrium state under the on-ramp based policies could be derived.
\par
  The analysis of this study was a generalization of \citet{Fu2022-nl} in terms of considering the heterogeneity of the value of the schedule delay. We proved that the QRP holds under the condition of the schedule delay cost function, even if commuters have heterogeneity, and showed the applications of the QRP to welfare analysis. These facts indicate that the QRP is a powerful tool for analyzing/characterizing the departure time choice problem in corridor networks, with or without commuters' heterogeneity. This presentation of the robustness of the analysis methodology using the QRP is one of the contributions of this study.
\par
  One of the most important issues to address in future studies is to analyze whether the QRP holds, considering more general heterogeneity of commuters, for example, the values of time and preferred arrival time.
  Particularly when the heterogeneity of values of time is considered; this may lead to different results regarding Pareto improvement.
  Moreover, it is important to generalize the network structures.
  For the generalization of network structures, it is considered adequate to start with analyzing one-to-many corridor networks (i.e., evening commute problem) and tree networks without the route choice, which simple applications of this study can analyze.
  In fact, considering homogenous commuters, \citet{Fu2022-nl} showed that the QRP also holds in evening commute problems under a slightly different assumption of the schedule delay cost function to morning commute problems.
  Given this, it is natural to conjecture that positive results can be derived from heterogeneous commuters' cases.
\par
  Future work should also examine the relationship between the DUE and DSO states when the QRP does not hold. 
  The numerical experiments in \citet{Fu2022-nl}, assuming homogeneous commuters, showed that a DUE state exists even when the QRP does not hold.
  Similarly, a DUE state may exist in cases of heterogeneous commuters.
  However, when the QRP does not hold, this study's approach can not derive the DUE state. Therefore, developing analytical methods for such situations is a significant future work.

\section*{CRediT authorship contribution statement}
\textbf{Takara Sakai}: Methodology, Formal analysis, Writing-original draft, and Funding acquisition.
\textbf{Takashi Akamatsu}: Conceptualization, Methodology, Writing-review and editing, Supervision, Project administration, and Funding acquisition.
\textbf{Koki Satsukawa}: Methodology, Formal analysis, Writing-original draft, Writing-review and editing, Project administration, Supervision, and Funding acquisition.

\section*{Acknowledgements}
This work was supported by JSPS KAKENHI Grant Numbers JP20J21744, JP21H01448, and JP20H02267.

\newpage

\appendix
\renewcommand{\thethm}{\Alph{section}.\arabic{thm}}
\renewcommand{\thedfn}{\Alph{section}.\arabic{dfn}}
\renewcommand{\thelem}{\Alph{section}.\arabic{lem}}
\renewcommand{\theasm}{\Alph{section}.\arabic{asm}}
\renewcommand{\thepro}{\Alph{section}.\arabic{pro}}
\renewcommand{\thecnj}{\Alph{section}.\arabic{cnj}}

\section{List of notations}\label{sec:ListofNotations}

\begin{table}[!ht]
  \centering
  \label{tab:ListofNotations}
  \begin{tabularx}{\linewidth}{lX}
      \hline
       \textit{Parameters}  & {}
       \\
       $N$ & Number of origin nodes/bottlenecks in the corridor network
       \\
       $\ClN = \{ i \mid i=1,2,...,N \}$
       & Set of origin nodes/bottlenecks
       \\
       $d_{i}$ & Free flow travel time from origin $i$ to the destination
       \\
       $\mu_{i}$ & Capacity of bottleneck $i$
       \\ 
       $\overline{\mu}_{i} = \mu_{i} - \mu_{(i+1)}$ &
       Difference between the capacities of bottleneck i and the upstream bottleneck $i + 1$, where $\mu_{(N+1)} = 0$
       \\
       $K$ & Number of commuter groups with different value of schedule delay
       \\
       $\ClK = \{ k \mid k=1,2,...,K \}$ & Set of commuter groups
       \\
       $Q_{i,k}$ & Total mass of $(i,k)$-commuters
       \\
       $\ClT$ & Set of arrival times
       \\
       $t^{\mathrm{P}} = 0$ & Preferred arrival time of all commuters
       \\
       $c(t)$ & Base schedule delay cost function
       \\
       $\beta^{k} \in (0, 1]$ & Parameter representing the heterogeneity of commuters with $\beta^{K}<... <\beta^{1}=1$
       \\
       $\beta^{k}c(t)$ & Schedule delay cost function of the group $k$ commuters as shown in \cref{fig:NetworkAndSDF}(b)
       \\
       $\overline{\beta}^{k} = \beta^{k} - \beta^{(k+1)}$ & Difference between the parameters representing the heterogeneity of commuters, where $\overline{\beta}^{(K+1)} = \beta^{K}$
       \\
       $\overline{c}(T)$ & 
       Schedule delay cost as a function of the arrival time window $T$ (\cref{fig:InverseSDF})
       \\
       $\Gamma(T)$ & 
       Set of arrival times related to $\overline{c}(T)$ defined by \eqref{eq:Gamma} (\cref{fig:InverseSDF})
       \\
       {} & {}
       \\
       \textit{Variables}
       \\
       $A_{i}(t)$ & Cumulative arrival flows for bottleneck $i$ by time $t$
       \\
       $D_{i}(t)$ & Cumulative departure flows for bottleneck $i$ by time $t$
       \\
       $\tau_{i}(t)$ &
       Arrival time at bottleneck $i$ for commuters whose destination arrival time is $t$
       \\
       $\sigma_{i}(t)$ &
        Departure time at bottleneck $i$ for commuters whose destination arrival time is $t$
       \\
       $w_{i}(t)$ & Queuing delay at bottleneck $i$ for the commuters whose destination arrival time is $t$
       \\
       $\rho_{i,k}$ & Equilibrium commuting cost of $(i,k)$-commuters
       \\
       $p_{i}(t)$ & Optimal price pattern at bottleneck $i$
       \\
       $x(\sigma_{i}(t))$
       & Departure flow rate at link $i$ for commuters whose destination arrival time is $t$
       \\
       $q_{i,k}(t)$ &
       Inflow rate to the network of $(i, k)$-commuters whose destination arrival time is $t$ 
       \\
       {} & {}
       \\
       \textit{Superscripts meaning of variable $X$} 
       \\
       $\S{X}$ & $X$ in the dynamic system optimal state (DSO)
       \\
       $\E{X}$ & $X$ in the dynamic user equilibrium state (DUE)
       \\
       $\PBP{X}$ & $X$ in the equilibrium under the partial bottleneck pricing (PBP)
       \\
       $\RM{X}$ & $X$ in the equilibrium under the on-ramp metering (RM)
       \\
       $\RP{X}$ & $X$ in the equilibrium under the on-ramp pricing (RP)
       \\
       $\PRM{X}$ &  $X$ in the equilibrium under the partial on-ramp metering (PRM)
       \\
       $\PRP{X}$ & $X$ in the equilibrium under the partial on-ramp pricing (PRP)
       \\
       \hline
    \end{tabularx}%
  \end{table}%

\section{Proofs of propositions}\label{sec:Proofs}

\subsection{Proof of \cref{lem:S_i_subset_S_i+1}}\label{subsec:Prf_Pro_Lem_S_i_subset_S_i+1}
We first prove the following two lemmas:
\begin{lem}
	Suppose that \cref{asm:p>0q>0} holds,
	there exists time $t\in\mathcal{T}$ when the following condition holds for all $k \in \ClK$ and for all $i \in \ClN \setminus \{ N \}$
	\begin{align}
		\S{q}_{i, k}(t) > 0 \quad \mbox{and} \quad
		\S{q}_{(i+1), k}(t) > 0.  
	\end{align}
	\label{lem:DSO_q>0q>0}
\end{lem}
\begin{prf}[\cref{lem:DSO_q>0q>0}]
	Proof by contradiction.
  Consider arbitrary $i \in \ClN \setminus \{ N \}$ and $k \in \ClK$.
  Suppose that the following relationships hold:
		\begin{align}
			&\S{q}_{i, k}(t) > 0 \quad \Rightarrow \quad \S{q}_{(i+1), k}(t) = 0
			&&\forall t \in \ClT,
			\label{eq:Lem_DSO_q>0q>0_1}
			\\
			&\S{q}_{(i+1), k}(t) > 0 \quad \Rightarrow \quad \S{q}_{i, k}(t) = 0
			&&\forall t \in \ClT.
			\label{eq:Lem_DSO_q>0q>0_2}
		\end{align}
			Let $t' \in \ClT$ be the time such that $\S{q}_{i, k}(t') > 0$.
			At time $t'$, the following conditions hold by \cref{eq:Lem_DSO_q>0q>0_1} and the optimality condition \eqref{eq:DSO-OC-DTC}:
		\begin{align} 
				&\S{\rho}_{i, k} = \beta^{k} c(t') + \sum_{j;j\leq i} \S{p}_{j}(t') + d_{i},
				\\
				&\S{\rho}_{(i+1), k} \leq \beta^{k} c(t') + \sum_{j;j\leq (i+1)} \S{p}_{j}(t') + d_{(i+1)}
				&&\because 
				\mbox{condition \eqref{eq:Lem_DSO_q>0q>0_1}}.
		\end{align}
		By combining the two conditions, we obtain
		\begin{align}
			&\S{\rho}_{(i+1), k} \leq \S{\rho}_{i, k} + \S{p}_{(i+1)}(t') + d_{(i+1)} - d_{i}
			\quad \Rightarrow \quad 
			\S{\rho}_{(i+1), k} \leq \S{\rho}_{i, k} 
			&&\because \S{p}_{(i+1)}(t') \geq 0, \quad d_{(i+1)} - d_{i} > 0.
			\label{eq:Lem_DSO_q>0q>0_Contradiction_1}
		\end{align}
		\par
		In the same way, 	for arbitrary $i \in \ClN \setminus \{ N \}$ and $k \in \ClK$, let $t'' \in \ClT $ be the time such that $\S{q}_{(i+1), k}(t'') > 0$.
		At time $t''$, the following conditions hold by \cref{eq:Lem_DSO_q>0q>0_2} and the optimality condition \eqref{eq:DSO-OC-DTC}:
		\begin{align}
				&\S{\rho}_{(i+1), k} = \beta^{k} c(t'') + \sum_{j;j\leq (i+1)} \S{p}_{j}(t'') + d_{(i+1)},
				\\
				&\S{\rho}_{i, k} \leq \beta^{k} c(t'') + \sum_{j;j\leq i} \S{p}_{j}(t'') + d_{i}
				&&\because \mbox{condition \eqref{eq:Lem_DSO_q>0q>0_2}}.
		\end{align}
	  Therefore, we have
		\begin{align}
			&\S{\rho}_{i, k} \leq \S{\rho}_{(i+1), k} - \S{p}_{(i+1)}(t'') + d_{i} - d_{(i+1)}
			\quad \Rightarrow \quad
			\S{\rho}_{i, k} < \S{\rho}_{(i+1), k}
      \notag
      \\
			&\hspace{65mm}
      \because \S{p}_{(i+1)}(t'') > 0  \mbox{ (\cref{asm:p>0q>0})}, \quad d_{(i+1)} - d_{i} > 0.
			\label{eq:Lem_DSO_q>0q>0_Contradiction_2}
		\end{align}
		\par
		\cref{eq:Lem_DSO_q>0q>0_Contradiction_1,eq:Lem_DSO_q>0q>0_Contradiction_2} contradict each other.
		This completes the proof of \cref{lem:DSO_q>0q>0}.\qed
\end{prf}

\begin{lem}
	If \cref{asm:p>0q>0} holds, then
	the equilibrium commuting cost satisfies the following inequality:
	\begin{align}
		&\S{\rho}_{i, k'} -d_{i} 
    < 
    \S{\rho}_{(i+1), k'} -d_{(i+1)}.
	\end{align}
	\label{lem:DSO_rho<rho}
\end{lem}
\begin{prf}[\cref{lem:DSO_rho<rho}]
	Consider time $t'$ such that $\S{q}_{i, k'}(t') > 0$ and $\S{q}_{(i+1), k'}(t) > 0$, 
	which always exists by \cref{lem:DSO_q>0q>0}.
	The equilibrium commuting costs of  $(i, k')$-commuters and $(i+1, k')$-commuters are represented as 
	\begin{align}
		&\S{\rho}_{i, k'} 
		= \beta^{k'} c(t') + \sum_{j;j \leq i} \S{p}_{j}(t') + d_{i},
		\\
		&\S{\rho}_{(i+1), k'} 
		= \beta^{k'} c(t') + \sum_{j;j \leq (i+1)} \S{p}_{j}(t')
    + d_{(i+1)}.
	\end{align}
	Then, 
	\begin{align}
		&\S{\rho}_{(i+1), k'} - \S{\rho}_{i, k'} 
    - d_{(i+1)} + d_{i}
    = \S{p}_{(i+1)}(t')
		> 0
		&&\because \mbox{\cref{asm:p>0q>0} and optimality conditions of [DSO]}.
	\end{align}
	This completes the proof. \qed
\end{prf}
	We can now prove \cref{lem:S_i_subset_S_i+1}.
	Condition \eqref{eq:Lem_S_i_subset_S_i+1} is equivalent to the following condition:
	\begin{align}
		&\S{q}_{i}(t) > 0
		\quad \Rightarrow \quad
		\S{q}_{(i+1)}(t) > 0
		&&\forall i \in \ClN \setminus \{ N \}.
	\end{align}
	Considering origin $i \in \ClN \setminus \{ N \}$, 
	suppose that there exists time $t'$ when 
	the following condition holds:
	\begin{align}
		\exists t' \in \ClT_{i}
		\quad \S{q}_{i}(t')>0 \quad \mbox{and} \quad \S{q}_{(i+1)}(t') = 0.
	\end{align}
	From the optimality condition, the following formula holds:
	\begin{align}
		\exists k' \in \ClK
		\quad \S{\rho}_{i, k'}(t') = c^{k'}(t') + \sum_{j;j\leq i} \S{p}_{j}(t') + d_{i}.
	\end{align}
		Consider $(i+1, k')$-commuters.
		Because $\S{q}_{(i+1)}(t') = 0$, we find $\S{q}_{(i+1), k'}(t') = 0$.
		Thus, we have
	\begin{align}
		&c^{k'}(t') + \sum_{j;j \leq (i+1)} \S{p}_{j}(t') 
    + d_{(i+1)}
		\geq \S{\rho}_{(i+1), k'}
		&&\because \mbox{optimality condition}
		\\
		\Leftrightarrow \quad
		&c^{k'}(t') + \sum_{j;j \leq i} \S{p}_{j}(t')
    + d_{(i+1)}
		\geq \S{\rho}_{(i+1), k'}
		&&\because \mbox{\cref{asm:p>0q>0}}
		\\
		\Leftrightarrow \quad
		&\S{\rho}_{i, k'} - d_{i}
    \geq \S{\rho}_{i+1, k'} - d_{(i+1)}.
	\end{align}
		This contradicts $\S{\rho}_{i, k'} -d_{i} 
    < 
    \S{\rho}_{(i+1), k'} -d_{(i+1)} $
		 (\cref{lem:DSO_rho<rho}).
		Then, \cref{lem:S_i_subset_S_i+1} is proved.
		\qed

\subsection{Proof of \cref{lem:q=mu}}\label{subsec:Prf_Pro_Lem_q=mu}
	We prove \cref{eq:Lem_q=mu} by induction.
	It follows from \cref{asm:p>0q>0} and 
	the optimality conditions that
	\begin{align}
		\S{q}_{N}(t) > 0
		\quad \Leftrightarrow \quad
		\S{p}_{N}(t) > 0
		\quad \Rightarrow \quad
		\S{q}_{N}(t) = \mu_{N} (=\overline{\mu}_{N}).
	\end{align}
	Thus, \cref{eq:Lem_q=mu} holds for $i = N$.
	\par
	Subsequently, we assume that \cref{eq:Lem_q=mu} holds for 
	$i=i'$:
	\begin{align}
		\S{q}_{i'}(t) 
		&=
		\begin{dcases}
			\overline{\mu}_{i'}
			\quad &\mathrm{if}\quad t \in \ClT_{i'}
			\\
			0 \quad &\mathrm{other}
		\end{dcases}
	\end{align}
Considering $i = i' - 1$, 
it follows from \cref{asm:p>0q>0} and 
the optimality conditions that
\begin{align}
	\S{q}_{(i'-1)}(t) > 0
	\quad \Leftrightarrow \quad
	&\S{p}_{(i'-1)}(t) > 0
	&&\because \mbox{\cref{asm:p>0q>0}}
	\\
	\quad \Rightarrow \quad
	&\S{q}_{(i'-1)}(t) = \mu_{(i'-1)} - \sum_{j;j > (i'-1)} \S{q}_{j}(t)
	&&\because \mbox{optimality conditions of [DSO]}
	\\
	\quad \Rightarrow \quad
	&\S{q}_{(i'-1)}(t) = \mu_{(i'-1)} - \mu_{i'} (= \overline{\mu}_{i'}) 
	&&\because \mbox{\cref{lem:S_i_subset_S_i+1}}
\end{align}
Thus, \cref{eq:Lem_q=mu} holds for $i = i' - 1$. 
\par
  Using the demand conservation condition \eqref{eq:DSO_ODcncv}, 
  we obtain the following time window length $\S{T}_{i} = | \S{\ClT}_{i} |$:
  \begin{align}
    &\S{T}_{i} = \dfrac{\sum_{k \in \ClK} Q_{i,k}}{\mu_{(i+1)} - \mu_{i}}
		&&\forall i \in \ClN.
  \end{align}
This completes the proof. \qed

\subsection{Proof of \cref{lem:LateSubmodular}}\label{subsec:Prf_lem:LateSubmodular}
	Condition \eqref{eq:LateSubmodular} can be rewritten as 
	\begin{align}
		\left(\beta^{k} - \beta^{(k+1)} \right)
		\left( 
			c(t) - c(t') 
		\right)< 0.
	\end{align}
	This holds because $\beta^{k} - \beta^{(k+1)}>0$, and $c(\cdot)$ is an increasing function of $t \geq 0$.
	This completes the proof of \cref{lem:LateSubmodular}.
	\qed

\subsection{Proof of \cref{lem:EarlySupermodular}}
\label{subsec:Prf_lem:EarlySupermodular}
	Condition \eqref{eq:EarlySupermodular} can be rewritten as 
	\begin{align}
		\left(\beta^{k} - \beta^{(k+1)} \right)
		\left( 
			c(t) - c(t') 
		\right)> 0.
	\end{align}
	This holds because $\beta^{k} - \beta^{(k+1)}<0$ and $c(\cdot)$ is a decreasing function of $t<0$.
	This completes the proof of \cref{lem:EarlySupermodular}.
	\qed

\subsection{Proofs of \cref{thm:QRP_Corridor} and \cref{pro:DUE_Solution}}
\label{subsec:Prf_thm:QRP_Corridor}
	We simultaneously prove \cref{thm:QRP_Corridor} and \cref{pro:DUE_Solution} by deriving the analytical solution to the DUE problem under \cref{eq:QRP_Corridor}.
Specifically, assuming \cref{eq:QRP_Corridor},
we derive the equilibrium commuting cost $\E{\rho}_{i, k}$ and arrival rate $\E{q}_{i, k}(t)$ as we guarantee no contradictions with DUE conditions.
\par
	First, we prove the following lemmas, which are useful for proving \cref{thm:QRP_Corridor} and \cref{pro:DUE_Solution}.
\begin{lem}
  Consider $i \in \ClN$ and $t^{-}, t^{+} \in \S{\ClT}_{i}$, $t^{-} \neq t^{+}$.
  If $c(t^{-}) = c(t^{+})$, then, there exists $k \in \ClK$ such that $t^{-}, t^{+} \in \S{\ClT}_{i, k}$.
  \label{lem:c=c_SamaGroup}
\end{lem}
\begin{prf}[\cref{lem:c=c_SamaGroup}]
    Because $t^{-} [t^{+}] \in \S{\ClT}_{i}$, there exists $k^{-} [k^{+}]  \in \ClK$ such that $t^{-} \in \ClT_{i,k^{-}} [t^{+} \in \ClT_{i,k^{+}}]$.
    Then, considering the DSO solution and the optimality condition \eqref{eq:DSO-OC-DTC} of the DSO problem, we have
  \begin{align}
    &k^{-} = \argmin_{k \in \ClK}. \left\{ \beta^{k}c(t^{-}) + \sum_{j;j\leq i} \S{p}_{i}(t^{-}) + d_{i}\right\}
    \quad \mbox{and} \quad 
    k^{+} = \argmin_{k \in \ClK}. \left\{ \beta^{k}c(t^{+}) + \sum_{j;j\leq i} \S{p}_{i}(t^{+}) + d_{i}\right\}.
  \end{align}
    Hence, if $c(t^{-}) = c(t^{+})$, it must be $k^{-} = k^{+}$. 
    This completes the proof.
    \qed
\end{prf}
\begin{lem}
  Consider $i \in \ClN$ and $t^{-}, t^{+} \in \S{\ClT}_{i}$, $t^{-} \neq t^{+}$.The base schedule delay cost function and optimal pricing pattern in the DSO state satisfy the following relationship:
	\begin{align}
		&c(t^{-}) = c(t^{+})
		\quad \Rightarrow \quad 
		\S{p}_{i}(t^{-}) = \S{p}_{i}(t^{+})
		&&\forall i \in \ClN, \quad \forall t^{-}, t^{+} \in \ClT, \quad t^{-} \neq t^{+}.
    \label{eq:c=c_p=p}
	\end{align}
	\label{lem:c=c_p=p}
\end{lem}
\begin{prf}[\cref{lem:c=c_p=p}]
  According to \cref{lem:c=c_SamaGroup}, there exists $k \in \ClK$ such that $t^{-}, t^{+} \in \S{\ClT}_{i, k}$.
    Then, considering the optimality condition \eqref{eq:DSO-OC-DTC} of the DSO problem, we obtain 
 \begin{align}
    &\beta^{k}c(t^{-}) + \sum_{j;j\leq i} \S{p}_{i}(t^{-})  + d_{i} = \S{\rho}_{i,k}
    \quad \mbox{and} \quad
    \beta^{k}c(t^{+}) + \sum_{j;j\leq i} \S{p}_{i}(t^{+}) 
    + d_{i} = \S{\rho}_{i,k}.
  \end{align}
	Hence, 
	\begin{align}
		&c(t^{-}) = c(t^{+})
		\quad \Rightarrow \quad
		\sum_{j; j \leq i} \S{p}_{i}(t^{-}) =  \sum_{j; j \leq i} \S{p}_{i}(t^{+})
		&&\forall t \in \ClT_{i}.
		\label{eq:Htr_c=c_sump=sump}
	\end{align}
	Because \cref{eq:Htr_c=c_sump=sump} holds for all $i \in \ClN$, we can derive the relationship \eqref{eq:c=c_p=p} by induction.
	This completes the proof.
  \qed
\end{prf}
\par
Now, \cref{thm:QRP_Corridor} and \cref{pro:DUE_Solution} can be proved.
We first show that $\E{w}_{i}(t) (= \S{p}_{i}(t))$ satisfies the consistency condition \cref{eq:DUE_queue_consistency}.
\begin{align}
	&\sum_{j;j \leq i} \E{\dot{w}}(t) = 
	\sum_{j;j \leq i} \E{\dot{p}}(t) 
	= 
	\begin{dcases}
		- \beta^{k}\dot{c}(t) < 1
		\quad &\mathrm{if}\quad t \in \S{\ClT}_{i}
		\\
		0
		\quad &\mathrm{if}\quad  t \notin \S{\ClT}_{i}
	\end{dcases}
	&&\forall i \in \ClN.
\end{align}
Thus, the consistency condition holds because of \cref{eq:QRP_Corridor}.
\par
Second, we show that $\E{w}_{i}(t) (= \S{p}_{i}(t))$ can satisfy 
the departure time condition \cref{eq:DUE-DTC}, 
compared with the optimality condition \cref{eq:DSO-OC-DTC}:
\begin{align}
	&\begin{dcases}
		\sum_{j; j\leq i} \E{w}_{j}(t) + \beta^{k} c(t) 
    + d_{i} = \E{\rho}_{i, k} 
		\quad &\mathrm{if}\quad \E{q}_{i, k}(t) > 0
		\\
		\sum_{j; j\leq i} \E{w}_{j}(t) + \beta^{k} c(t) 
    + d_{i} \geq \E{\rho}_{i, k}
		\quad &\mathrm{if}\quad \E{q}_{i, k}(t) = 0
	 \end{dcases}
	 &&\forall k \in \ClK,
	 \quad \forall i \in \mathcal{N},
	 \quad \forall t \in \ClT,
	 \\
	 \Leftrightarrow \quad 
	&\begin{dcases}
		\sum_{j; j\leq i} \S{p}_{j}(t) + \beta^{k} c(t)  + d_{i} = \E{\rho}_{i, k} 
		\quad &\mathrm{if}\quad \E{q}_{i, k}(t) > 0
		\\
		\sum_{j; j\leq i} \S{p}_{j}(t) + \beta^{k} c(t)  + d_{i} \geq \E{\rho}_{i, k}
		\quad &\mathrm{if}\quad \E{q}_{i, k}(t) = 0
	 \end{dcases}
    &&\forall k \in \ClK,
    \quad \forall i \in \mathcal{N},
    \quad \forall t \in \ClT.
\end{align}
Therefore
\begin{align}
	&\E{\rho}_{i, k} = \S{\rho}_{i, k}
	&&\forall k \in \ClK, \quad \forall i \in \ClN.
\end{align}
Third, we show that the arrival rate $\E{q}_{i, k}(t)$ satisfies 
the non-negativity condition.
The flow pattern in the DUE state $\E{q}_{i}(t)$ 
can be represented in the following form:
\begin{align}
	\E{q}_{i}(t) = 
	\sum_{j; j \geq i} \E{q}_{j}(t) - \sum_{j;j > i} \E{q}_{j}(t)
	&=\mu_{i} \E{\dot{\sigma}}_{i}(t) - 
	\mu_{(i+1)} \E{\dot{\sigma}}_{(i+1)}(t)
	\\
	&=
    \begin{dcases}
      \mu_{i} - \mu_{(i+1)} \E{\dot{\sigma}}_{(i+1)}(t)
      \quad &\mathrm{if}\quad t \notin \E{\ClT}_{(i-1)}
      \quad \mbox{[Case-1]}
      \\
    \mu_{i} \E{\dot{\sigma}}_{i}(t) - 
    \mu_{(i+1)} \E{\dot{\sigma}}_{(i+1)}(t)
    \quad &\mathrm{if}\quad t \in \E{\ClT}_{(i-1)}
    \quad \mbox{[Case-2]}
    \end{dcases}
\end{align}
Considering [Case-1], the nonnegative flow condition holds:
\begin{align}
	\mu_{i} - \mu_{(i+1)} \E{\dot{\sigma}_{(i+1)}}
  &\geq 
  \mu_{i} - \mu_{(i+1)} - \mu_{(i+1)}\beta^{\max} \dot{c}(t)
  \notag
  \\
  &> \mu_{i} - \mu_{(i+1)} - \mu_{(i+1)}\beta^{\max} 
  \dfrac{\mu_{i} - \mu_{(i+1)}}{ \beta^{\max} \mu_{(i+1)}}
	= 0
  &&\because \mbox{\cref{eq:QRP_Corridor}}.
\end{align}
\color{black}
Considering [Case-2], 
we define the function $f(\beta^{k'}, \beta^{k''}, t)$ as follows:
\begin{align}
	f(\beta^{k'}, \beta^{k''}, t)
	&\equiv 
	\mu_{i} \E{\dot{\sigma}}_{i}(t) - 
	\mu_{(i+1)} \E{\dot{\sigma}}_{(i+1)}(t) 
	= 
	\mu_{i} \left(1 + \beta^{k'}\dot{c}(t) \right) - 
	\mu_{(i+1)} \left(1 + \beta^{k''}\dot{c}(t) \right).
\end{align}
	Considering $t < 0 \Rightarrow \dot{c}(t) < 0$ and $t \geq 0 \Rightarrow \dot{c}(t) \geq 0$, the minimum value of the function $f(\beta^{k'}, \beta^{k''}, t)$ is calculated as follows:
  \begin{align}
    &
      \min_{\beta^{k'}, \beta^{k''}, t}.
      f(\beta^{k'}, \beta^{k''}, t)
      = 
        \min. \left\{ 
          \min_{t;t<0}. f(\beta^{\max}, \beta^{\min}, t),
          \quad
          \min_{t;t\geq 0}. f(\beta^{\min}, \beta^{\max}, t)
        \right\}.
  \end{align}
  Based on the above, we find that the minimum value of the function $f(\beta^{k'}, \beta^{k''}, t)$ is nonnegative as follows:
  \begin{align}
    f(\beta^{\max}, \beta^{\min}, t) 
    &= 
    \mu_{i} \left(1 + \beta^{\max} \dot{c}(t) \right) - 
    \mu_{(i+1)} \left(1 + \beta^{\min}\dot{c}(t) \right)
    \notag
    \\
    &=
    \mu_{i} - \mu_{(i+1)}
    + 
    \dot{c}(t) \left( \mu_{i} \beta^{\max} - \mu_{(i+1)}\beta^{\min} \right)
    \notag
    \\
    &>
    \mu_{i} - \mu_{(i+1)}
    + \max \left\{ 
     - 1, \quad 
     -
     \dfrac{\mu_{i} - \mu_{(i+1)}}{\mu_{i}\beta^{\max} - \mu_{(i+1)}\beta^{\min}}
    \right\}
    \left( \mu_{i} \beta^{\max} - \mu_{(i+1)}\beta^{\min} \right)
    &&\because \mbox{\cref{eq:QRP_Corridor}}
    \notag
    \\
    &\geq 0
    &&\forall t < 0.
  \end{align}
  \begin{align}
    f(\beta^{\min}, \beta^{\max}, t) 
    &= 
    \mu_{i} \left(1 + \beta^{\min} \dot{c}(t) \right) - 
    \mu_{(i+1)} \left(1 + \beta^{\max}\dot{c}(t) \right)
    \notag
    \\
    &=
    \mu_{i} - \mu_{(i+1)}
    + 
    \dot{c}(t) \left( \mu_{i} \beta^{\min} - \mu_{(i+1)}\beta^{\max} \right)
    \notag
    \\
    &>
    \mu_{i} - \mu_{(i+1)} + 
    \dfrac{\mu_{i} - \mu_{(i+1)}}{ \beta^{\max} \mu_{(i+1)}}
    \left( \mu_{i} \beta^{\min} - \mu_{(i+1)}\beta^{\max} \right)
    &&\because \mbox{\cref{eq:QRP_Corridor}}
    \notag
    \\
    &= 
    \mu_{i} - \mu_{(i+1)} + \left( \mu_{i} - \mu_{(i+1)} \right)\dfrac{\mu_{i} \beta^{\min}}{\mu_{(i+1)}\beta^{\max}}
    - \left( \mu_{i} - \mu_{(i+1)} \right)
    \notag
    \\
    &=\left( \mu_{i} - \mu_{(i+1)} \right)\dfrac{\mu_{i} \beta^{\min}}{\mu_{(i+1)}\beta^{\max}}
    > 0
    &&\forall t \geq 0.
  \end{align}
  \color{black}
\noindent
	Thus, the nonnegative flow condition holds.
\par
	Finally, we confirm that the demand conservation 
	condition holds by \cref{lem:c=c_p=p}:
\begin{align}
	\int_{\ClT_{i, k}} \E{q}_{i, k}(t) \mathrm{d} t 
	&= 
	\int_{\ClT_{i, k}} 
	\left(
		\mu_{i} \E{\dot{\sigma}}_{i}(t) 
	- \mu_{(i+1)} \E{\dot{\sigma}}_{(i+1)}(t)
	\right) \mathrm{d} t
	\notag
	\\
	&=
	\int_{\ClT_{i, k}} 
	\left(
		\mu_{i}     \left( 1 - \sum_{j;j<i}     \dot{\E{w}}_{j}(t) \right) 
	- \mu_{(i+1)} \left( 1 - \sum_{j;j<(i+1)} \dot{\E{w}}_{j}(t) \right) 
	\right) \mathrm{d} t
	\notag
	\\
	&= 
	\int_{\ClT_{i, k}}
	\left(
		\mu_{i} - \mu_{(i+1)}
	\right) \mathrm{d} t
	- \mu_{i}
	\int_{\ClT_{i, k}}
		\sum_{j;j<i} \E{\dot{w}}_{j}(t)
		\mathrm{d} t
	+
	\mu_{(i+1)}
	\int_{\ClT_{i, k}}
	\sum_{j;j<(i+1)} \E{\dot{w}}_{j}(t) \mathrm{d} t
	\notag
	\\
	&= 
	\sum_{l;l\leq k} Q_{i, l}
	- \mu_{i}
	\int_{\ClT_{i, k}}
		\sum_{j;j<i} \S{\dot{p}}_{j}(t)
		\mathrm{d} t
	+
	\mu_{(i+1)}
	\int_{\ClT_{i, k}}
	\sum_{j;j<(i+1)} \S{\dot{p}}_{j}(t) \mathrm{d} t
	\notag
	\\
	&= \sum_{l;l\leq k} Q_{i, l}
	\qquad \qquad \qquad \forall k \in \ClK, \quad \forall i \in \ClN, 
	\quad \because\mbox{\cref{lem:c=c_p=p}}.
\end{align}
This completes the proof.
\qed

		
		

	

\subsection{Proof of \cref{lem:PBP_Solution}}\label{subsec:Prf_lem:PBP_Solution}
  We show that the variable set in \cref{eq:PBP_q,eq:PBP_rho,eq:PBP_w} satisfies the equilibrium condition of [DUE-PBP].
  First, we show that $\{ \PBP{\Vtw}(t) \}_{t \in \ClT}$ in \cref{eq:PBP_w} satisfies the consistency condition \eqref{eq:DUE-PBP_ConsistencyCondition}:
\begin{align}
  1 - 
  \sum_{j; j \in \overline{\ClN}^{\mathrm{P}}, j \leq i} 
  \PBP{\dot{w}}_{j}(t) 
  &= 
  \begin{dcases}
    1 - 
    \sum_{j; j \leq i} 
    \E{\dot{w}}_{j}(t) > 0
    \quad &\mathrm{if}\quad i \in \overline{\ClN}^{\mathrm{P}}
    \\
    1 - 
    \sum_{j; j < \PBP{i}} 
    \E{\dot{w}}_{j}(t) > 0
    \quad &\mathrm{if}\quad i \in \ClN^{\mathrm{P}}
  \end{dcases}
\end{align}
where $\PBP{i} = \argmin. \ClN^{\mathrm{P}}$
\par
  Subsequently, we show that $\{ \PBP{\Vtq}(t) \}_{t \in \ClT}$ in \cref{eq:PBP_q} guarantees the nonnegative flow condition as follows:
\begin{align}
    \PBP{q}_{i}(t) &= 
    \mu_{i} \PBP{\dot{\sigma}}_{i}(t) 
    -  
    \mu_{(i+1)} \PBP{\dot{\sigma}}_{(i+1)}(t)
    \\
    &=
    \mu_{i}
    \left(1 -  \sum_{j; j < \min. \{i, \PBP{i} \}} \E{\dot{w}}_{j}(t) \right)
    -
    \mu_{(i+1)}
    \left(1 -  \sum_{j; j \leq \min. \{i, \PBP{i} \}}  \E{\dot{w}}_{j}(t) \right)
    \\
    &=
    \begin{dcases}
      \left(  \mu_{i} - \mu_{(i+1)} \right)
      \left( 1 - \sum_{j; j \leq \PBP{i} }  \E{\dot{w}}_{j}(t) \right) > 0
      \quad &\mathrm{if}\quad i > \PBP{i}
      \\
      \E{q}_{i}(t) > 0
      \quad &\mathrm{if}\quad i \leq \PBP{i}
    \end{dcases}
\end{align}
\par
  Because the commuting cost of each commuter equals that in the DSO and DUE states, we confirm that $\{ \PBP{\Vtq}(t) \}_{t \in \ClT}$, $\{ \PBP{\Vtw}(t) \}_{t \in \ClT}$, and $\{ \PBP{\Vtrho} \}$ satisfy the departure time condition \eqref{eq:DUE-PBP_DTC}.
  Moreover, the queueing condition \eqref{eq:DUE-PBP_Queue} holds because the queuing delay pattern at the no-pricing bottleneck $i \in \ClN \setminus \ClN^{\mathrm{P}}$ equals that in the DUE state.
\par
Finally, the demand conservation condition holds as follows:
		\begin{align}
			\int_{\ClT_{i, k}}
			\PBP{q}_{i}(t)
			\mathrm{d} t
			&= \int_{\ClT_{i, k}}
			\left(
				\mu_{i} \PBP{\dot{\sigma}}_{i}(t) 
				-  
				\mu_{(i+1)} \PBP{\dot{\sigma}}_{(i+1)}(t)
			\right)
			\mathrm{d} t
			\\
			&=
			\int_{\ClT_{i, k}}
			\left(
			\mu_{i}
			\left(1 -  \sum_{j; j < \min. \{i, \PBP{i} \}} \E{\dot{w}}_{j}(t) \right)
			-
			\mu_{(i+1)}
			\left(1 -  \sum_{j; j \leq \min. \{i, \PBP{i} \}}  \E{\dot{w}}_{j}(t) \right)
			\right)
			\mathrm{d} t
			\\
			&=
			\begin{dcases}
        \int_{\ClT_{i, k}}
				\left(  \mu_{i} - \mu_{(i+1)} \right)
				\left( 1 - \sum_{j; j \leq \PBP{i} }  \E{\dot{w}}_{j}(t) \right)  \mathrm{d} t
				\quad &\mathrm{if}\quad i > \PBP{i}
				\\
        \int_{\ClT_{i, k}}
				\E{q}_{i}(t)
        \mathrm{d} t
				\quad &\mathrm{if}\quad i \leq \PBP{i}
			\end{dcases}
			\\
			&=\sum_{l; l \leq k} Q_{i, l}
			\qquad \qquad  \qquad  \qquad \qquad \qquad  
      \forall k \in \ClK, \quad \forall i \in \ClN.
		\end{align}
		This completes the proof.
		\qed

		\subsection{Proof of \cref{thm:PBP_Pareto}}\label{subsec:Prf_thm:PBP_Pareto}
		Comparing the analytical solutions
		(DSO: \cref{pro:DSO_Solution}, 
		DUE: \cref{pro:DUE_Solution},
		and PBP: \cref{lem:PBP_Solution}), 
		we obtain \cref{eq:PBP_EquilibriumCost_are_equal}.
		\par
		According to \cref{lem:PBP_Solution}, we have
		\begin{align}
			\PBP{Z}(\ClN^{\mathrm{P}}) &= 
			\sum_{i \in \ClN} \sum_{k \in \ClK}  Q_{i,k} \E{\rho}_{i,k} 
			- \sum_{i \in \ClN^{\mathrm{P}}} 
			\int_{\ClT} \mu_{i} \S{p}_{i}(t) \mathrm{d} t
      \\
			&= 
			\sum_{i \in \ClN} \sum_{k \in \ClK}  Q_{i,k} \E{\rho}_{i,k} 
			- \sum_{i \in \ClN} \int_{\ClT} \mu_{i} \S{p}_{i}(t) \mathrm{d} t
			+ \sum_{i \in \overline{\ClN}^{\mathrm{P}}} 
			\int_{\ClT} \mu_{i} \S{p}_{i}(t) \mathrm{d} t.
		\end{align}
		Therefore, inequality \eqref{eq:PBP_TotalCost_DSO<PBP} holds, and the equality in \cref{eq:PBP_TotalCost_DSO<PBP} holds if and only if $\ClN^{\mathrm{P}}=\ClN$.s
		\qed
		\subsection{Proof of \cref{lem:RM_solution}}\label{subsec:prf_lem:RM_solution}
			We show that $\{ \RM{q}_{i^{\ast}, k}(t) \}$, 
			$\{ \RM{\rho}_{i^{\ast}, k} \}$ and $\{ \RM{w}_{i^{\ast}}(t) \}$ in \cref{lem:RM_solution} satisfy the equilibrium condition of [DUE-RM].
		First, the demand conservation condition holds as follows:
		\begin{align}
			\int_{\ClT} \RM{q}_{i^{\ast}, k}(t) \mathrm{d}t 
			&=
			\int_{ \E{\ClT}_{i, k} \setminus \E{\ClT}_{i, (k-1)}} 
			\overline{\mu}_{i}  \mathrm{d}t 
			= Q_{i,k}
			&&\forall k \in \ClK, \quad \forall i^{\ast} \in \ClN^{\ast}.
		\end{align}
		The nonnegative condition of flows also holds.
		Because $\RM{w}_{i^{\ast}}(t) = \sum_{j;j\leq i} \E{w}_{j}(t)$, 
		the nonnegative condition of queues and the consistency condition hold.
		Moreover, the departure time choice condition also holds because the queueing cost of $(i,k)$-commuters arriving at time $t$ at the destination equals that in the DUE state.
		Finally, the queueing condition holds as follows:
		\begin{align}
			&\RM{w}_{i^{\ast}}(t) > 0
			\quad \Rightarrow \quad
			t \in \E{\ClT}_{i}
			\quad \Rightarrow \quad
			\sum_{k \in \ClK} \RM{q}_{i^{\ast}, k}(t) = \overline{\mu}_{i},
			\\
			&\RM{w}_{i^{\ast}}(t) = 0
			\quad \Rightarrow \quad
			t \notin \E{\ClT}_{i}
			\quad \Rightarrow \quad
			\sum_{k \in \ClK} \RM{q}_{i^{\ast}, k}(t) = 0 < \overline{\mu}_{i}.
		\end{align}
		Thus 	$\{ \RM{q}_{i^{\ast}, k}(t) \}$, 
		$\{ \RM{\rho}_{i^{\ast}, k} \}$,
		$\{ \RM{w}_{i}(t) \}$,
		$\{ \RM{w}_{i^{\ast}}(t) \}$
		are solutions to [DUE-RM].
		This completes the proof.
		\qed

	\subsection{Proof of \cref{lem:RP_solution}}\label{subsec:prf_lem:RP_solution}
		We show that 
			$\{ \RP{q}_{i^{\ast}, k}(t) \}$, 
			$\{ \RP{\rho}_{i^{\ast}, k} \}$, and 
			$\{ \RP{w}_{i^{\ast}}(t) \}$
			in \cref{lem:RP_solution}
			satisfy the equilibrium conditions of [DUE-RP].
		First, the demand conservation condition holds as follows:
		\begin{align}
			\int_{\ClT} \RP{q}_{i^{\ast}, k}(t) \mathrm{d}t 
			&=
			\int_{ \E{\ClT}_{i, k} \setminus \E{\ClT}_{i, (k-1)}} 
			\overline{\mu}_{i}  \mathrm{d}t 
			= Q_{i,k}
			&&\forall k \in \ClK, \quad \forall i^{\ast} \in \ClN^{\ast}.
		\end{align}
		The nonnegative condition of flows also holds.
		Because $\RP{w}_{i^{\ast}}(t) = 0$, 
		the nonnegative condition of queues and consistency condition hold.
		Moreover, 
		the departure time choice condition also holds because the queueing cost of $(i,k)$-commuters arriving at time $t$ at the destination equals that in the DUE state.
		Finally, the queueing condition holds as follows:
		\begin{align}
			&\RP{w}_{i}(t) = 0
			\quad \Rightarrow \quad
			t \notin \E{\ClT}_{i}
			\quad \Rightarrow \quad
			\sum_{k \in \ClK} \RP{q}_{i^{\ast}, k}(t) \leq \overline{\mu}_{i}.
		\end{align}
		Thus 	
		$\{ \RP{q}_{i^{\ast}, k}(t) \}$, 
		$\{ \RP{\rho}_{i^{\ast}, k} \}$,
		$\{ \RP{w}_{i^{\ast}}(t) \}$
		are solutions to [DUE-RM].
		This completes the proof. \qed

	\subsection{Proofs of \cref{lem:PRM_Solution} and \cref{lem:PRP_Solution}}\label{subsec:prf_PartialPolicies}
		We can prove \cref{lem:PRM_Solution,lem:PRP_Solution} by an approach similar to that of the proof of \cref{lem:PBP_Solution}.
		Specifically, by substituting the variables shown in each lemma into each equilibrium conditions, we can confirm that the variables are solutions to [DUE-PRM] and [DUE-PRP], respectively.
		\qed

\subsection{Proof of \cref{thm:TotalCost_DSO=RP<RM=DUE}}
\label{subsec:Prf_thm:TotalCost_DSO=RP<RM=DUE}
Comparing the analytical solutions
(DSO: \cref{pro:DSO_Solution}, 
DUE: \cref{pro:DUE_Solution},
RM: \cref{lem:RM_solution}, and RP: \cref{lem:RP_solution}), 
we obtain \cref{eq:EquilibriumCost_are_equal_Full}.
\par
Using \cref{pro:DSO_Solution}, \cref{pro:DUE_Solution},
\cref{lem:RM_solution}, and \cref{lem:RP_solution}, $\S{Z}$, $\E{Z}$, $\RM{Z}$, and $\RP{Z}$ are calculated as follows:
\begin{align}
	&\S{Z}= 
	\sum_{i \in \ClN} \sum_{k \in \ClK} Q_{i,k} \S{\rho}_{i,k} 
	- \sum_{i \in \ClN} \int_{\ClT} \mu_{i} \S{p}_{i}(t) \mathrm{d} t
	= 
	\sum_{i \in \ClN} \sum_{k \in \ClK}  Q_{i,k} \E{\rho}_{i,k} 
	- \sum_{i \in \ClN} \int_{\ClT} \mu_{i} \S{p}_{i}(t) \mathrm{d} t,
	\\
	&\E{Z}
	=\sum_{i \in \ClN} \sum_{k \in \ClK} Q_{i,k} \E{\rho}_{i,k},
	\\
	&\RM{Z} = 
	\sum_{i \in \ClN} \sum_{k \in \ClK}  Q_{i,k} \RM{\rho}_{i,k} 
	= 
	\sum_{i \in \ClN} \sum_{k \in \ClK}  Q_{i,k} \E{\rho}_{i,k},
	\\
	&\RP{Z} = 
	\sum_{i \in \ClN} \sum_{k \in \ClK} Q_{i,k} \RP{\rho}_{i,k} 
	- \sum_{i \in \ClN} \int_{\ClT} \overline{\mu}_{i} \S{P}_{i}(t) \mathrm{d} t
	= 
	\sum_{i \in \ClN} Q_{i} \S{\rho}_{i}
	- 
	\sum_{i \in \ClN} \int_{\ClT} \mu_{i} \S{p}_{i}(t) \mathrm{d} t.
\end{align}
Thus, we have inequality \eqref{eq:S-BP-RM-E_CostIneq}.\qed

\subsection{Proof of \cref{thm:TotalCost_DSO<PRP<PRM=DUE}}
\label{subsec:Prf_thm:TotalCost_DSO<PRP<PRM=DUE}
Comparing the analytical solutions
(DSO: \cref{pro:DSO_Solution}, 
DUE: \cref{pro:DUE_Solution},
PRM: \cref{lem:PRM_Solution},
and PRP: \cref{lem:PRP_Solution}), 
we obtain \cref{eq:EquilibriumCost_are_equal_Partial_1}.
\par
From \cref{lem:PRM_Solution,lem:PRP_Solution}, we obtain
	\begin{align}
		&\PRM{Z}(\ClN^{\ast\mathrm{P}}) 
		= \sum_{i \in \ClN} \sum_{k \in \ClK} \PRM{\rho}_{i,k} Q_{i,k}
		= \E{Z},
		\\
		&\PRP{Z}(\ClN^{\ast\mathrm{P}}) 
		= \sum_{i \in \ClN} \sum_{k \in \ClK} \PRP{\rho}_{i,k} Q_{i,k}
		= \S{Z} + \sum_{i \in\ClN \setminus \ClN^{\ast\mathrm{P}}} \int_{\ClT} w_{i}(t) \mu_{i} \mathrm{d}t.
	\end{align}
	Thus, inequalities \eqref{eq:TotalCost_DSO<PRP<DUE} and \eqref{eq:TotalCost_DSO<PRM=DUE} hold.
    Note that $\S{Z} = \PRP{Z}(\ClN^{\ast\mathrm{P}})$ holds if 
    and only if $\ClN^{\mathrm{P}}=\ClN$, and 
    $\PRP{Z}(\ClN^{\ast\mathrm{P}}) = \E{Z}$ holds if 
    and only if $\ClN^{\mathrm{P}}=\emptyset$.
	\qed


\begin{thebibliography}{40}
\expandafter\ifx\csname natexlab\endcsname\relax\def\natexlab#1{#1}\fi
\providecommand{\url}[1]{\texttt{#1}}
\providecommand{\href}[2]{#2}
\providecommand{\path}[1]{#1}
\providecommand{\DOIprefix}{doi:}
\providecommand{\ArXivprefix}{arXiv:}
\providecommand{\URLprefix}{URL: }
\providecommand{\Pubmedprefix}{pmid:}
\providecommand{\doi}[1]{\href{http://dx.doi.org/#1}{\path{#1}}}
\providecommand{\Pubmed}[1]{\href{pmid:#1}{\path{#1}}}
\providecommand{\bibinfo}[2]{#2}
\ifx\xfnm\relax \def\xfnm[#1]{\unskip,\space#1}\fi
\bibitem[{Akamatsu and Wada(2017)}]{Akamatsu2017-bi}
\bibinfo{author}{Akamatsu, T.}, \bibinfo{author}{Wada, K.},
  \bibinfo{year}{2017}.
\newblock \bibinfo{title}{Tradable network permits: A new scheme for the most
  efficient use of network capacity}.
\newblock \bibinfo{journal}{Transportation Research Part C: Emerging
  Technologies} \bibinfo{volume}{79}, \bibinfo{pages}{178--195}.
\bibitem[{Akamatsu et~al.(2015)Akamatsu, Wada and Hayashi}]{Akamatsu2015-ip}
\bibinfo{author}{Akamatsu, T.}, \bibinfo{author}{Wada, K.},
  \bibinfo{author}{Hayashi, S.}, \bibinfo{year}{2015}.
\newblock \bibinfo{title}{The corridor problem with discrete multiple
  bottlenecks}.
\newblock \bibinfo{journal}{Transportation Research Part B: Methodological}
  \bibinfo{volume}{81}, \bibinfo{pages}{808--829}.
\bibitem[{Akamatsu et~al.(2021)Akamatsu, Wada, Iryo and
  Hayashi}]{Akamatsu2021-zg}
\bibinfo{author}{Akamatsu, T.}, \bibinfo{author}{Wada, K.},
  \bibinfo{author}{Iryo, T.}, \bibinfo{author}{Hayashi, S.},
  \bibinfo{year}{2021}.
\newblock \bibinfo{title}{A new look at departure time choice equilibrium
  models with heterogeneous users}.
\newblock \bibinfo{journal}{Transportation Research Part B: Methodological}
  \bibinfo{volume}{148}, \bibinfo{pages}{152--182}.
\bibitem[{Arnott(1998)}]{Arnott1998-of}
\bibinfo{author}{Arnott, R.}, \bibinfo{year}{1998}.
\newblock \bibinfo{title}{Congestion tolling and urban spatial structure}.
\newblock \bibinfo{journal}{Journal of regional science} \bibinfo{volume}{38},
  \bibinfo{pages}{495--504}.
\bibitem[{Arnott and DePalma(2011)}]{Arnott2011-rb}
\bibinfo{author}{Arnott, R.}, \bibinfo{author}{DePalma, E.},
  \bibinfo{year}{2011}.
\newblock \bibinfo{title}{The corridor problem: Preliminary results on the
  no-toll equilibrium}.
\newblock \bibinfo{journal}{Transportation Research Part B: Methodological}
  \bibinfo{volume}{45}, \bibinfo{pages}{743--768}.
\bibitem[{Arnott et~al.(1988)Arnott, de~Palma and Lindsey}]{Arnott1988-pl}
\bibinfo{author}{Arnott, R.}, \bibinfo{author}{de~Palma, A.},
  \bibinfo{author}{Lindsey, R.}, \bibinfo{year}{1988}.
\newblock \bibinfo{title}{Schedule delay and departure time decisions with
  heterogeneous commuters}.
\newblock \bibinfo{journal}{Transportation research record}
  \bibinfo{volume}{476}, \bibinfo{pages}{56--57}.
\bibitem[{Arnott et~al.(1990)Arnott, de~Palma and Lindsey}]{Arnott1990-ta}
\bibinfo{author}{Arnott, R.}, \bibinfo{author}{de~Palma, A.},
  \bibinfo{author}{Lindsey, R.}, \bibinfo{year}{1990}.
\newblock \bibinfo{title}{Departure time and route choice for the morning
  commute}.
\newblock \bibinfo{journal}{Transportation Research Part B: Methodological}
  \bibinfo{volume}{24}, \bibinfo{pages}{209--228}.
\bibitem[{Arnott et~al.(1992)Arnott, de~Palma and Lindsey}]{Arnott1992-jg}
\bibinfo{author}{Arnott, R.}, \bibinfo{author}{de~Palma, A.},
  \bibinfo{author}{Lindsey, R.}, \bibinfo{year}{1992}.
\newblock \bibinfo{title}{Route choice with heterogeneous drivers and
  group-specific congestion costs}.
\newblock \bibinfo{journal}{Regional Science and Urban Economics}
  \bibinfo{volume}{22}, \bibinfo{pages}{71--102}.
\bibitem[{Arnott et~al.(1993)Arnott, de~Palma and Lindsey}]{Arnott1993-mq}
\bibinfo{author}{Arnott, R.}, \bibinfo{author}{de~Palma, A.},
  \bibinfo{author}{Lindsey, R.}, \bibinfo{year}{1993}.
\newblock \bibinfo{title}{Properties of dynamic traffic equilibrium involving
  bottlenecks, including a paradox and metering}.
\newblock \bibinfo{journal}{Transportation Science} \bibinfo{volume}{27},
  \bibinfo{pages}{148--160}.
\bibitem[{Arnott et~al.(1994)Arnott, de~Palma and Lindsey}]{Arnott1994-dy}
\bibinfo{author}{Arnott, R.}, \bibinfo{author}{de~Palma, A.},
  \bibinfo{author}{Lindsey, R.}, \bibinfo{year}{1994}.
\newblock \bibinfo{title}{The welfare effects of congestion tolls with
  heterogeneous commuters}.
\newblock \bibinfo{journal}{Journal of Transport Economics and Policy}
  \bibinfo{volume}{28}, \bibinfo{pages}{139--161}.
\bibitem[{van~den Berg and Verhoef(2011)}]{Van_den_Berg2011-qb}
\bibinfo{author}{van~den Berg, V.}, \bibinfo{author}{Verhoef, E.T.},
  \bibinfo{year}{2011}.
\newblock \bibinfo{title}{Congestion tolling in the bottleneck model with
  heterogeneous values of time}.
\newblock \bibinfo{journal}{Transportation Research Part B: Methodological}
  \bibinfo{volume}{45}, \bibinfo{pages}{60--78}.
\bibitem[{van~den Berg and Verhoef(2010)}]{Van_den_Berg2010-ty}
\bibinfo{author}{van~den Berg, V.A.C.}, \bibinfo{author}{Verhoef, E.T.},
  \bibinfo{year}{2010}.
\newblock \bibinfo{title}{Why congestion tolling could be good for the
  consumer: The effects of heterogeneity in the values of schedule delay and
  time on the effects of tolling}.
\newblock \bibinfo{journal}{Tinbergen Institute Discussion Paper 2010-016/3} .
\bibitem[{Chen et~al.(2015)Chen, Nie and Yin}]{Chen2015-ku}
\bibinfo{author}{Chen, H.}, \bibinfo{author}{Nie, Y.m.}, \bibinfo{author}{Yin,
  Y.}, \bibinfo{year}{2015}.
\newblock \bibinfo{title}{Optimal {Multi-Step} toll design under general user
  heterogeneity}.
\newblock \bibinfo{journal}{Transportation Research Procedia}
  \bibinfo{volume}{7}, \bibinfo{pages}{341--361}.
\bibitem[{Daganzo(1985)}]{Daganzo1985-ls}
\bibinfo{author}{Daganzo, C.F.}, \bibinfo{year}{1985}.
\newblock \bibinfo{title}{The uniqueness of a time-dependent equilibrium
  distribution of arrivals at a single bottleneck}.
\newblock \bibinfo{journal}{Transportation Science} \bibinfo{volume}{19},
  \bibinfo{pages}{29--37}.
\bibitem[{Daniel et~al.(2009)Daniel, Gisches and Rapoport}]{Daniel2009-la}
\bibinfo{author}{Daniel, T.E.}, \bibinfo{author}{Gisches, E.J.},
  \bibinfo{author}{Rapoport, A.}, \bibinfo{year}{2009}.
\newblock \bibinfo{title}{Departure times in {Y-Shaped} traffic networks with
  multiple bottlenecks}.
\newblock \bibinfo{journal}{American Economic Review} \bibinfo{volume}{99},
  \bibinfo{pages}{2149--2176}.
\bibitem[{DePalma and Arnott(2012)}]{DePalma2012-la}
\bibinfo{author}{DePalma, E.}, \bibinfo{author}{Arnott, R.},
  \bibinfo{year}{2012}.
\newblock \bibinfo{title}{Morning commute in a single-entry traffic corridor
  with no late arrivals}.
\newblock \bibinfo{journal}{Transportation Research Part B: Methodological}
  \bibinfo{volume}{46}, \bibinfo{pages}{1--29}.
\bibitem[{Fu et~al.(2022)Fu, Akamatsu, Satsukawa and Wada}]{Fu2022-nl}
\bibinfo{author}{Fu, H.}, \bibinfo{author}{Akamatsu, T.},
  \bibinfo{author}{Satsukawa, K.}, \bibinfo{author}{Wada, K.},
  \bibinfo{year}{2022}.
\newblock \bibinfo{title}{Dynamic traffic assignment in a corridor network:
  Optimum versus equilibrium}.
\newblock \bibinfo{journal}{Transportation Research Part B: Methodological}
  \bibinfo{volume}{161}, \bibinfo{pages}{218--246}.
\bibitem[{Hall(2018)}]{Hall2018-lg}
\bibinfo{author}{Hall, J.D.}, \bibinfo{year}{2018}.
\newblock \bibinfo{title}{Pareto improvements from lexus lanes: The effects of
  pricing a portion of the lanes on congested highways}.
\newblock \bibinfo{journal}{Journal of Public Economics} .
\bibitem[{Hall(2021)}]{Hall2021-xo}
\bibinfo{author}{Hall, J.D.}, \bibinfo{year}{2021}.
\newblock \bibinfo{title}{Can tolling help everyone? estimating the aggregate
  and distributional consequences of congestion pricing}.
\newblock \bibinfo{journal}{Journal of the European Economic Association}
  \bibinfo{volume}{19}, \bibinfo{pages}{441--474}.
\bibitem[{Hendrickson and Kocur(1981)}]{Hendrickson1981-cu}
\bibinfo{author}{Hendrickson, C.}, \bibinfo{author}{Kocur, G.},
  \bibinfo{year}{1981}.
\newblock \bibinfo{title}{Schedule delay and departure time decisions in a
  deterministic model}.
\newblock \bibinfo{journal}{Transportation Science} \bibinfo{volume}{15},
  \bibinfo{pages}{62--77}.
\bibitem[{Iryo and Yoshii(2007)}]{Iryo2007-ne}
\bibinfo{author}{Iryo, T.}, \bibinfo{author}{Yoshii, T.}, \bibinfo{year}{2007}.
\newblock \bibinfo{title}{Equivalent optimization problem for finding
  equilibrium in the bottleneck model with departure time choices}, in:
  \bibinfo{booktitle}{4th {IMA} International Conference on Mathematics in
  {TransportInstitute} of Mathematics and its Applications},
  \bibinfo{publisher}{trid.trb.org}.
\bibitem[{Kuwahara(1990)}]{Kuwahara1990-nk}
\bibinfo{author}{Kuwahara, M.}, \bibinfo{year}{1990}.
\newblock \bibinfo{title}{Equilibrium queueing patterns at a {Two-Tandem}
  bottleneck during the morning peak}.
\newblock \bibinfo{journal}{Transportation Science} \bibinfo{volume}{24},
  \bibinfo{pages}{217--229}.
\bibitem[{Kuwahara and Akamatsu(1993)}]{Kuwahara1993-vt}
\bibinfo{author}{Kuwahara, M.}, \bibinfo{author}{Akamatsu, T.},
  \bibinfo{year}{1993}.
\newblock \bibinfo{title}{Dynamic equilibrium assignment with queues for a
  one-to-many {OD} pattern}, \bibinfo{publisher}{In: Daganzo, C.F. (Ed.),
  Proceedings of the 12th International Symposium on Transportation and Traffic
  Theory. Elsevior, Berkeley}. pp. \bibinfo{pages}{185--204}.
\bibitem[{Lago and Daganzo(2007)}]{Lago2007-wq}
\bibinfo{author}{Lago, A.}, \bibinfo{author}{Daganzo, C.F.},
  \bibinfo{year}{2007}.
\newblock \bibinfo{title}{Spillovers, merging traffic and the morning commute}.
\newblock \bibinfo{journal}{Transportation Research Part B: Methodological}
  \bibinfo{volume}{41}, \bibinfo{pages}{670--683}.
\bibitem[{Laih(1994)}]{Laih1994-hi}
\bibinfo{author}{Laih, C.H.}, \bibinfo{year}{1994}.
\newblock \bibinfo{title}{Queueing at a bottleneck with single- and multi-step
  tolls}.
\newblock \bibinfo{journal}{Transportation Research Part A: Policy and
  Practice} \bibinfo{volume}{28}, \bibinfo{pages}{197--208}.
\bibitem[{Li and Huang(2017)}]{Li2017-xg}
\bibinfo{author}{Li, C.Y.}, \bibinfo{author}{Huang, H.J.},
  \bibinfo{year}{2017}.
\newblock \bibinfo{title}{Morning commute in a single-entry traffic corridor
  with early and late arrivals}.
\newblock \bibinfo{journal}{Transportation Research Part B: Methodological}
  \bibinfo{volume}{97}, \bibinfo{pages}{23--49}.
\bibitem[{Li et~al.(2020)Li, Huang and Yang}]{Li2020-zz}
\bibinfo{author}{Li, Z.C.}, \bibinfo{author}{Huang, H.J.},
  \bibinfo{author}{Yang, H.}, \bibinfo{year}{2020}.
\newblock \bibinfo{title}{Fifty years of the bottleneck model: A bibliometric
  review and future research directions}.
\newblock \bibinfo{journal}{Transportation Research Part B: Methodological}
  \bibinfo{volume}{139}, \bibinfo{pages}{311--342}.
\bibitem[{Lindsey(2004)}]{Lindsey2004-aw}
\bibinfo{author}{Lindsey, R.}, \bibinfo{year}{2004}.
\newblock \bibinfo{title}{Existence, uniqueness, and trip cost function
  properties of user equilibrium in the bottleneck model with multiple user
  classes}.
\newblock \bibinfo{journal}{Transportation Science} \bibinfo{volume}{38},
  \bibinfo{pages}{293--314}.
\bibitem[{Lindsey et~al.(2012)Lindsey, van~den Berg and
  Verhoef}]{Lindsey2012-gg}
\bibinfo{author}{Lindsey, R.}, \bibinfo{author}{van~den Berg, V.A.C.},
  \bibinfo{author}{Verhoef, E.T.}, \bibinfo{year}{2012}.
\newblock \bibinfo{title}{Step tolling with bottleneck queuing congestion}.
\newblock \bibinfo{journal}{Journal of Urban Economics} \bibinfo{volume}{72},
  \bibinfo{pages}{46--59}.
\bibitem[{Liu et~al.(2015)Liu, Nie and Hall}]{Liu2015-ay}
\bibinfo{author}{Liu, Y.}, \bibinfo{author}{Nie, Y.m.}, \bibinfo{author}{Hall,
  J.}, \bibinfo{year}{2015}.
\newblock \bibinfo{title}{A semi-analytical approach for solving the bottleneck
  model with general user heterogeneity}.
\newblock \bibinfo{journal}{Transportation Research Part B: Methodological}
  \bibinfo{volume}{71}, \bibinfo{pages}{56--70}.
\bibitem[{Luenberger(1997)}]{Luenberger1997-la}
\bibinfo{author}{Luenberger, D.G.}, \bibinfo{year}{1997}.
\newblock \bibinfo{title}{Optimization by Vector Space Methods}.
\newblock \bibinfo{publisher}{John Wiley {\&} Sons}.
\bibitem[{Newell(1987)}]{Newell1987-wk}
\bibinfo{author}{Newell, G.F.}, \bibinfo{year}{1987}.
\newblock \bibinfo{title}{The morning commute for nonidentical travelers}.
\newblock \bibinfo{journal}{Transportation Science} \bibinfo{volume}{21},
  \bibinfo{pages}{74--88}.
\bibitem[{Osawa et~al.(2018)Osawa, Fu and Akamatsu}]{Osawa2018-hg}
\bibinfo{author}{Osawa, M.}, \bibinfo{author}{Fu, H.},
  \bibinfo{author}{Akamatsu, T.}, \bibinfo{year}{2018}.
\newblock \bibinfo{title}{First-best dynamic assignment of commuters with
  endogenous heterogeneities in a corridor network}.
\newblock \bibinfo{journal}{Transportation Research Part B: Methodological}
  \bibinfo{volume}{117}, \bibinfo{pages}{811--831}.
\bibitem[{Papageorgiou and Kotsialos(2002)}]{Papageorgiou2002-jj}
\bibinfo{author}{Papageorgiou, M.}, \bibinfo{author}{Kotsialos, A.},
  \bibinfo{year}{2002}.
\newblock \bibinfo{title}{Freeway ramp metering: an overview}.
\newblock \bibinfo{journal}{IEEE Transactions on Intelligent Transportation
  Systems} \bibinfo{volume}{3}, \bibinfo{pages}{271--281}.
\bibitem[{Rachev and R{\"{u}}schendorf(1998)}]{Rachev1998-bb}
\bibinfo{author}{Rachev, S.T.}, \bibinfo{author}{R{\"{u}}schendorf, L.},
  \bibinfo{year}{1998}.
\newblock \bibinfo{title}{Mass Transportation Problems: Volume I: Theory}.
\newblock \bibinfo{publisher}{Springer Science {\&} Business Media}.
\bibitem[{Sakai et~al.(2022)Sakai, Satsukawa and Akamatsu}]{Sakai2022-vm}
\bibinfo{author}{Sakai, T.}, \bibinfo{author}{Satsukawa, K.},
  \bibinfo{author}{Akamatsu, T.}, \bibinfo{year}{2022}.
\newblock \bibinfo{title}{Non-existence of queues for system optimal departure
  patterns in tree networks}.
\newblock \bibinfo{journal}{arXiv preprint}
  \href{http://arxiv.org/abs/2205.06015}{\tt arXiv:2205.06015}.
\bibitem[{Smith(1984)}]{Smith1984-wr}
\bibinfo{author}{Smith, M.J.}, \bibinfo{year}{1984}.
\newblock \bibinfo{title}{The existence of a {Time-Dependent} equilibrium
  distribution of arrivals at a single bottleneck}.
\newblock \bibinfo{journal}{Transportation Science} \bibinfo{volume}{18},
  \bibinfo{pages}{385--394}.
\bibitem[{Takayama and Kuwahara(2017)}]{Takayama2017-gx}
\bibinfo{author}{Takayama, Y.}, \bibinfo{author}{Kuwahara, M.},
  \bibinfo{year}{2017}.
\newblock \bibinfo{title}{Bottleneck congestion and residential location of
  heterogeneous commuters}.
\newblock \bibinfo{journal}{Journal of Urban Economics} \bibinfo{volume}{100},
  \bibinfo{pages}{65--79}.
\bibitem[{Vickrey(1969)}]{Vickrey1969-ic}
\bibinfo{author}{Vickrey, W.S.}, \bibinfo{year}{1969}.
\newblock \bibinfo{title}{Congestion theory and transport investment}.
\newblock \bibinfo{journal}{American Economic Review} \bibinfo{volume}{59},
  \bibinfo{pages}{251--260}.
\bibitem[{Wada and Akamatsu(2013)}]{Wada2013-li}
\bibinfo{author}{Wada, K.}, \bibinfo{author}{Akamatsu, T.},
  \bibinfo{year}{2013}.
\newblock \bibinfo{title}{A hybrid implementation mechanism of tradable network
  permits system which obviates path enumeration: An auction mechanism with
  day-to-day capacity control}.
\newblock \bibinfo{journal}{Transportation Research Part E: Logistics and
  Transportation Review} \bibinfo{volume}{60}, \bibinfo{pages}{94--112}.

\end{thebibliography}
\end{document}